\documentclass[a4paper,12pt,reqno]{amsart}

\usepackage[T1]{fontenc}
\usepackage[utf8]{inputenc}
\usepackage{lmodern}
\usepackage[english]{babel}

\usepackage{microtype} 
\usepackage[%
	left=2.5cm,       
	right=2.5cm,      
	top=3.5cm,        
	bottom=3.5cm,     
	heightrounded,    
	bindingoffset=0mm 
]{geometry}

\usepackage{amsmath}
\usepackage{amssymb}
\usepackage{mathtools}
\usepackage{mathrsfs}

\usepackage[
colorlinks=true,
urlcolor=teal,
citecolor=magenta,
linkcolor=teal,
breaklinks=true]
{hyperref}

\usepackage[noabbrev,capitalize,nameinlink]{cleveref}

\usepackage{bookmark}

\usepackage[initials,msc-links,
]{amsrefs}

\usepackage[dvipsnames]{xcolor}
\usepackage{graphicx}
\usepackage{url}

\usepackage{enumitem}

\numberwithin{equation}{section}

\theoremstyle{plain}
\newtheorem{theorem}{Theorem}[section]
\newtheorem{lemma}[theorem]{Lemma}
\newtheorem{proposition}[theorem]{Proposition}
\newtheorem{corollary}[theorem]{Corollary}

\theoremstyle{definition}
\newtheorem{definition}[theorem]{Definition}
\newtheorem{remark}[theorem]{Remark}

\newcommand{\atoms}{\mathcal{X}}
\newcommand{\bdata}{\mathcal B_\Omega}
\newcommand{\di}{\,\mathrm{d}}
\newcommand{\e}{\varepsilon}
\newcommand{\mat}[4]{\begin{pmatrix}
#1&#2\\[.5ex]#3&#4
\end{pmatrix}}
\newcommand{\N}{\mathbb{N}}
\newcommand{\R}{\mathbb{R}}

\newcommand{\smat}[4]{\begin{psmallmatrix}
#1&#2\\[.5ex]#3&#4
\end{psmallmatrix}}
\newcommand{\sobs}{H^s_{\bar{u}}(\Omega;B_L)}
\newcommand{\sphere}{\mathbb{S}^1}
\newcommand{\svec}[2]{\begin{psmallmatrix}
#1 \\[.5ex] #2
\end{psmallmatrix}}
\newcommand{\xflute}[1]{\xrightarrow{\flute(#1)}}

\newcommand{\vect}[2]{\begin{pmatrix}
#1 \\[.5ex] #2
\end{pmatrix}}

\newcommand{\Z}{\mathbb{Z}}

\DeclareMathOperator{\cu}{curl}
\DeclareMathOperator{\diverg}{div}
\DeclareMathOperator{\GL}{GL}
\DeclareMathOperator{\flute}{flat}

\DeclareMathOperator{\jac}{Jac}
\DeclareMathOperator{\lip}{Lip}
\DeclareMathOperator{\supp}{supp}

\DeclarePairedDelimiter{\set}{\{}{\}}

\allowdisplaybreaks

\begin{document}

\title[Topological singularities arising from fractional-gradient energies]{Topological singularities arising from fractional-gradient energies}

\author[R.~Alicandro]{Roberto Alicandro}
\address[R.~Alicandro]{Dipartimento di Matematica e Applicazioni ``Renato Caccioppoli'', Università degli Studi di Napoli ``Federico II'',
via Cintia, Monte S.\ Angelo, 80126 Napoli (NA), Italy}
\email{roberto.alicandro@unina.it}

\author[A.~Braides]{Andrea Braides}
\address[A.~Braides]{Scuola Internazionale Superiore di Studi Avanzati (SISSA), via Bonomea 265, 34136 Trieste (TS), Italy}
\email{abraides@sissa.it}

\author[M.~Solci]{Margherita Solci}
\address[M.~Solci]{Dipartimento di Architettura, Design e Urbanistica, Piazza Duomo 6, 07041 Alghero (SS), Italy}
\email{marghe@uniss.it}

\author[G.~Stefani]{Giorgio Stefani}
\address[G.~Stefani]{Scuola Internazionale Superiore di Studi Avanzati (SISSA), via Bonomea 265, 34136 Trieste (TS), Italy}
\email{gstefani@sissa.it {\normalfont or} giorgio.stefani.math@gmail.com}

\date{\today}



\thanks{%
\textit{Acknowledgements}. 
The authors are members of INdAM and  GNAMPA.
R.\ Alicandro, A.\ Braides and M.\ Solci\ have received funding from INdAM under the INdAM--GNAMPA 2024 Project  
\textit{Asymptotic analysis of nonlocal variational problems} (grant agreement No.\ CUP\_E53\-C23\-001\-670\-001).
G.\ Stefani has received funding from INdAM under the INdAM--GNAMPA 2024 Project  \textit{Ottimizzazione e disuguaglianze funzionali per problemi geometrico-spettrali locali e non-locali} (grant agreement No.\ CUP\_E53\-C23\-001\-670\-001), and the INdAM--GNAMPA 2023 Project \textit{Problemi variazionali per funzionali e operatori non-locali} (grant agreement No.\ CUP\_E53\-C22\-001\-930\-001), and from the European Research Council (ERC) under the European Union’s Horizon 2020 research and innovation program (grant agreement No.~945655).%
}

\begin{abstract}
We prove that, on a planar regular domain, suitably scaled functionals of Ginzburg--Landau type, given by the sum of quadratic fractional Sobolev seminorms and a penalization term vanishing on the unitary sphere, $\Gamma$-converge to vortex-type energies with respect to the flat convergence of Jacobians. 
The compactness and the $\Gamma$-$\liminf$ follow by comparison with standard Ginzburg--Landau functionals depending on Riesz potentials.
The $\Gamma$-$\limsup$, instead, is achieved via a direct argument by joining a finite number of vortex-like functions suitably truncated around the singularity.
\end{abstract}

\maketitle

\section{Introduction}

\subsection{Classical framework}
\label{subseci:classical}

Let $\Omega\subset\R^2$ be a non-empty, connected, simply connected, bounded open set with Lipschitz boundary. 
Given a scale parameter $\e>0$ and an additional parameter $\lambda>0$, the \emph{Ginzburg--Landau functionals} 
\begin{equation*}
\GL_{\e,\lambda}(\,\cdot\,;\Omega)\colon H^1(\Omega;\R^2)\to[0,\infty]
\end{equation*}
are defined as 
\begin{equation}
\label{eqi:GL}
\GL_{\e,\lambda}(v;\Omega)
=
\frac{1}{|\log\e|}
\int_\Omega |Dv|^2\di x
+
\frac{\lambda}{\e^2|\log\e|}\int_\Omega\left(|v|^2-1\right)^2\di x
\end{equation}
for $v\in H^1(\Omega;\R^2)$.
One may prescribe a trace constraint at the boundary by imposing that $v|_{\partial\Omega}=g$ for some fixed boundary datum $g\colon\partial\Omega\to\sphere$ with (topological) degree $d=\deg(g,\partial\Omega)\in\Z$, in which case the functionals in~\eqref{eqi:GL} are restricted to the subspace
\begin{equation*}
H^1_g(\Omega;\R^2)
=
\set*{v\in H^1(\Omega;\R^2) : v|_{\partial\Omega}=g}.
\end{equation*}

Much effort has been devoted to understanding the asymptotic behavior of the minimizers~$v_\e$ of the functionals in~\eqref{eqi:GL} as the scale parameter vanishes. 
In the limit as $\e\to0^+$, the minimizers $v_\e$ of the functionals in~\eqref{eqi:GL} develop \emph{vortex-like singularities} of the form $\frac{x-x_i}{|x-x_i|}$ (possibly, up to a fixed rotation) for $|d|$ points $x_i$'s in $\overline\Omega$.

After the works~\cites{J99,JS02,S00,S98} (we also refer to the monograph~\cite{BBH94} and to~\cite{AP14} for a more detailed presentation of the problem, and to~\cite{ABO05} for the higher-dimensional setting), the picture is nowadays well understood.
The \emph{$\Gamma$-convergence} of the functionals in~\eqref{eqi:GL} as $\e\to0^+$ is related to the \emph{flat} (or \emph{$1$-Wasserstein}) \emph{convergence} of the Jacobians $\jac(v_\e)=\det(Dv_\e)$ of their minimizers $v_\e$ to an \emph{atomic} measure $\mu\in\atoms(\overline\Omega)$, where
\begin{equation}
\label{eqi:atoms}
\atoms(\overline\Omega)
=
\set*{
\sum_{i=1}^N d_i\delta_{x_i} : d_i\in\Z\ \text{and}\ x_i\in\overline\Omega\ \text{for}\ i=1,\dots, N,\ \text{with}\ N\in\N
}.
\end{equation}
In more precise terms, we can state the following result.
For a detailed presentation of the notion of \emph{$\Gamma$-convergence}, we refer to the monographs~\cites{B02,D93}.
For the definition of \emph{flat convergence}; i.e., in the dual norm with respect to Lipschitz functions, see \cref{subsec:flat}. 

\begin{theorem}[Compactness and $\Gamma$-convergence of Ginzburg--Landau energies]
\label{resi:GL}
Let $\Omega\subset\R^2$ be a non-empty, connected, simply connected, bounded open set with Lipschitz boundary, $g\in H^{\frac12}(\partial\Omega;\sphere)$ with $d=\deg(g|_{\partial\Omega},\partial\Omega)\in\Z$ and $\lambda>0$.

\begin{enumerate}[label=(\roman*),ref=(\roman*),itemsep=1ex,leftmargin=4ex,topsep=1ex]

\item[{\rm(i)}]
\emph{(Compactness)}
If $(v_{\e_k})_{k\in\N}\subset H^1_g(\Omega;\R^2)$, with $\e_k\to0^+$ as $k\to\infty$, is such that 
\begin{equation*}
\sup_{k\in\N}\GL_{\e_k,\lambda}(v_{\e_k};\Omega)<\infty,
\end{equation*}
then there exists a subsequence $(v_{\e_{k_j}})_{j\in\N}$ and $\mu\in\atoms(\overline\Omega)$ such that $\mu(\overline\Omega)=d$ and
\begin{equation*}
\jac(v_{\e_{k_j}})\,\mathscr L^2
\xflute{\overline\Omega}\mu
\quad
\text{as}\ j\to\infty.
\end{equation*}

\item[{\rm(ii)}]
\emph{($\Gamma$-$\liminf$ inequality)}
If $(v_{\e_k})_{k\in\N}\subset H^1_g(\Omega;\R^2)$, with $\e_k\to0^+$ as $k\to\infty$, is such that
\begin{equation}
\label{eqi:jacs_liminf}
\jac(v_{\e_k})\,\mathscr L^2
\xflute{\overline\Omega}\mu
\quad
\text{as}\ k\to\infty
\end{equation}
for some $\mu\in\atoms(\overline\Omega)$, then $\mu(\overline\Omega)=d$ and 
\begin{equation*}
\liminf_{k\to\infty}\GL_{\e_k,\lambda}(v_{\e_k};\Omega)
\ge 
2\pi|\mu|(\overline\Omega).
\end{equation*}

\item[{\rm(iii)}]
\emph{($\Gamma$-$\limsup$ inequality)}
If $\mu\in\atoms(\overline\Omega)$ is such that $\mu(\overline\Omega)=d$, then there exists a sequence $(v_{\e_k})_{k\in\N}\subset H^1_g(\Omega;\R^2)$, with $\e_k\to0^+$ as $k\to\infty$, such that~\eqref{eqi:jacs_liminf} holds and 
\begin{equation*}
\limsup_{k\to\infty}\GL_{\e_k,\lambda}(v_{\e_k};\Omega)
\le 
2\pi|\mu|(\overline\Omega).
\end{equation*}  

\end{enumerate}
\end{theorem}

\subsection{Fractional framework}

In the present work, we investigate the validity of a fractional analog of \cref{resi:GL}, in which the differential operator in the $L^2$ energy in~\eqref{eqi:GL} is replaced by the \emph{fractional} {(\emph{Riesz})} \emph{$s$-gradient}
\begin{equation}
\label{eqi:s-grad}
D^su(x)
=
(1-s)
\,
\frac{2^{s-1}}{\pi}\,\tfrac{\Gamma\left(\frac{3+s}{2}\right)}{\Gamma\left(\frac{3-s}{2}\right)}
\int_{\R^2}
\frac{(u(y)-u(x))\otimes(y-x)}{|y-x|^{3+s}}\di y,
\quad
x\in\R^2,
\end{equation}
for $s\in(0,1)$, where $\Gamma$ stands for \emph{Euler's Gamma function}.
The integro-differential operator in~\eqref{eqi:s-grad} has gained considerable interest in recent years, leading to a rapidly growing number of works concerning the development of a fractional analog of  classical calculus and of its applications.
For a non-comprehensive account on the existing literature, we refer to~\cites{SS15,SS18,CS19,CS23,CS22,BCCS22,BCM20,BCM21,KS22,S24,S20} and the references therein.

Our main aim is to replace~\eqref{eqi:GL} with functionals of the form
\begin{equation}
\label{eqi:e_s-GL}
\GL_{\e_s,\lambda}(u;\Omega)
=
\frac{1}{|\log\e_s|}
\int_\Omega |D^su|^2\di x
+
\frac{\lambda}{\e^2_s|\log\e_s|}\int_\Omega (|u|^2-1 )^2\di x
\end{equation}
for some scale parameter $\e_s>0$, depending on $s$, such that $\e_s\to0^+$ as $s\to1^-$, in such a way that the analog of Theorem \ref{resi:GL} holds.

In order to determine $\e_s$, we observe that the $s$-gradient in~\eqref{eqi:s-grad} can be equivalently presented as the gradient of the \emph{Riesz potential} of order $1-s$; that is, $D^s u
=
Dv$, where
\begin{equation}
\label{eqi:riesz}
v(x)
=
I_{1-s}u(x)
=
\frac{1-s}{1+s}
\,
\frac{2^{s-1}}{\pi}
\,
\tfrac{\Gamma\left(\frac{3+s}{2}\right)}{\Gamma\left(\frac{3-s}{2}\right)}
\int_{\R^2}\frac{u(y)}{|y-x|^{1+s}}\di y,
\quad
x\in\R^2.
\end{equation}
Therefore, since $I_{1-s}u$ tends to $ u$ as $s\to1^-$, we expect that, along a family of minimizers $u_s$, the functionals in~\eqref{eqi:e_s-GL} can be reasonably approximated as
\begin{equation}
\label{eqi:approx}
\GL_{\e_s,\lambda}(u_s;\Omega)
\sim
\frac{1}{|\log\e_s|}
\int_\Omega |Dv_s|^2\di x
+
\frac{\lambda}{\e^2_s|\log\e_s|}\int_\Omega (|v_s|^2-1 )^2\di x
\sim 
\GL_{\e_s,\lambda}(v_s;\Omega)
\end{equation}
with $v_s=I_{1-s}u_s$ a family of `almost minimizers' of the functionals in~\eqref{eqi:GL}.
Assuming that $u_s$ and $v_s$ are uniformly bounded and supported in a bounded neighborhood of $\Omega$, we expect that the size of the error in the approximation~\eqref{eqi:approx} should be not larger than
\begin{equation}
\label{eqi:approx2}
\frac{1}{\e^2_s|\log\e_s|}
\int_\Omega
\left|
(|v_s|^2-1)^2
-
(|u_s|^2-1)^2
\right|
\di x
\lesssim
\frac{(1-s)^2}{\e^2_s|\log\e_s|}
\,
[u_s]_{H^s}^2
\end{equation}
(we refer to \cref{res:L2_comp} in \cref{subsec:truncated_Riesz} below for the precise computations), where
\begin{equation}
\label{eqi:s-energy}
[u]_{H^s}^2
=
\int_{\R^2}\int_{\R^2}\frac{|u(y)-u(x)|^2}{|y-x|^{2+2s}}\di x\di y
\end{equation} 
is the \emph{Sobolev--Slobodeckij $H^s$-energy} of $u$.
Since the $L^2$ norm of the $s$-gradient in~\eqref{eqi:s-grad} is proportional to the energy in~\eqref{eqi:s-energy},
\begin{equation}
\label{eqi:s-equal}
\int_{\R^2}|D^s u|^2\di x
=
\frac{4^s}{2\pi}\,\frac{s\,\Gamma(1+s)}{\Gamma(2-s)}\,(1-s)\,[u]_{H^s}^2
\end{equation}
(e.g., see~\cite{SS15}*{Rem.~2.3}, as well as \cref{res:equiv-frac} below),
we get that the size of the error in the approximation~\eqref{eqi:approx} should be at most
\begin{equation*}
\frac{(1-s)}{\e^2_s|\log\e_s|}
\int_{\R^2}|D^s u_s|^2\di x
\sim
\frac{(1-s)}{\e^2_s|\log\e_s|}
\int_{\R^2}|Dv_s|^2\di x.
\end{equation*} 
Therefore, in order to re-absorb such an error, we have to require that its size is comparable to the one of the other term in $\GL_{\e_s,\lambda}(v_s)$; that is, 
\begin{equation*}
\frac{(1-s)}{\e^2_s|\log\e_s|}
\int_{\R^2}|Dv_s|^2\di x
\sim
\frac{1}{|\log\e_s|}
\int_{\R^2}|Dv_s|^2\di x,
\end{equation*}
from which we get that $\e_s\sim\sqrt{1-s}$.

\subsection{Statement of the main result}

Although rather naive, the approximation in~\eqref{eqi:approx} leads to the correct fractional analog of~\eqref{eqi:GL}. 
In view of~\eqref{eqi:s-energy} and of the identification in~\eqref{eqi:s-equal}, as customary we set 
\begin{equation*}
H^s(\R^2;\R^2)
=
\set*{u\in L^2(\R^2;\R^2)
:
[u]_{H^s}^2<\infty}
\end{equation*}
for $s\in(0,1)$.
Therefore, given $\lambda>0$, we define the \emph{Ginzburg--Landau fractional $s$-energies}
\begin{equation*}
\GL^s_\lambda(\,\cdot\,;\Omega)\colon
H^s(\R^2;\R^2)
\to
[0,\infty)
\end{equation*}   
by letting
\begin{equation}
\label{eqi:s-GL}
\GL^s_\lambda(u;\Omega)
=
\frac{1}{|\log(1-s)|}\int_{\mathbb R^2}|D^s u|^2\di x
+
\frac{\lambda}{(1-s)|\log(1-s)|}
\int_{\Omega}
(|u|^2-1 )^2
\di x
\end{equation}
for $u\in H^s(\R^2;\R^2)$.
By~\eqref{eqi:s-equal}, we can equivalently write 
\begin{equation*}
\begin{split}
\GL^s_\lambda(u;\Omega)
&=
\frac{1-s}{|\log(1-s)|}\frac2\pi
\,[u]^2_{H^s}
+
\frac{\lambda}{(1-s)|\log(1-s)|}
\int_{\Omega}
(|u|^2-1)^2
\di x
\end{split}
\end{equation*}
for $u\in H^s(\R^2;\R^2)$.
We note that analogous scalar energies (but with different scaling regimes) are related to a non-local approach to phase-transition problems, see~\cites{SV12,DFV18,DSV20}.

In order to state the fractional analog of \cref{resi:GL}, we need to introduce a suitable boundary condition, which---due to the non-locality of~\eqref{eqi:s-grad}---has to be prescribed on~$\R^2\setminus\Omega$ instead of just on $\partial\Omega$.
In addition, in view of the underlying approximation argument in~\eqref{eqi:approx}, we require that the boundary datum is bounded and supported in a bounded neighborhood of $\Omega$.
Hence, the set of boundary data we consider is defined as
\begin{equation}
\label{eqi:bdata}
\begin{split}
\bdata
=
\set*{
g\in H^1(\R^2;\R^2)\cap L^\infty(\R^2;\R^2)
:
\text{\parbox[c]{7cm}{\centering $g$ has compact support and \\ $|g|=1$ in an open neighborhood of $\partial\Omega$}}\,
}.
\end{split}
\end{equation}
Thus, given $g\in\bdata$ and $L\in [\|g\|_{L^\infty(\R^2)},\infty )$, for $s\in(0,1)$ we define the spaces
\begin{equation}
\label{eqi:sobs}
H^s_g(\Omega;B_L)
=
\set*{u\in H^s(\R^2;\R^2) : u=g\ \text{on}\ \R^2\setminus\Omega\ \text{and}\ \|u\|_{L^\infty}\le L}.
\end{equation}

In view of the representation of the $s$-gradient in~\eqref{eqi:s-grad} via the Riesz potential in~\eqref{eqi:riesz}, the approximation in~\eqref{eqi:approx}, and the analogy with \cref{resi:GL}, we are naturally led to define the {\em fractional Jacobian} of a function $u\in H^s(\R^2;\R^2)$ as the usual Jacobian of the function $v=I_{1-s}u$ for all $s\in(0,1)$, letting
\begin{equation}
\label{eqi:def_s-jac}
\jac^s(u)=\det(D^s u)=\det(Dv)=\jac(v).
\end{equation}
The convergence of $\jac^s(u)$ as $s\to1^-$ is then defined according to the customary definition of flat convergence.

With the above notation in force, our main result states as follows.

\begin{theorem}[Compactness and $\Gamma$-convergence of fractional Ginzburg--Landau energies] 
\label{resi:main}
Let $\Omega\subset\R^2$ be a non-empty, connected, simply connected, bounded open set with Lipschitz boundary, $g\in\bdata$ with $d=\deg(g|_{
\partial\Omega},\partial\Omega)\in\Z$, $L\in \left[\|g\|_{L^\infty},\infty\right)$ and $\lambda>0$.

\begin{enumerate}[label=(\roman*),ref=(\roman*),leftmargin=4ex,itemsep=1ex]

\item\label{item:main_comp} 
\emph{(Compactness)}
If $u_{s_k}\in H^{s_k}_g(\Omega;B_L)$, with $s_k\to1^-$ as $k\to\infty$, is such that 
\begin{equation*}
\GL^{s_k}_\lambda(u_{s_k};\Omega)<\infty,
\end{equation*} 
then there exists a subsequence $(u_{s_{k_j}})_{j\in\N}$ and $\mu\in\atoms(\overline\Omega)$ such that $\mu(\overline\Omega)=d$ and 
\begin{equation*}
\jac^{s_{k_j}}(u_{s_{k_j}})\xflute{\overline\Omega}\pi\mu
\quad
\text{as}\
j\to\infty.
\end{equation*}

\item\label{item:main_liminf} 
\emph{($\Gamma$-$\liminf$ inequality)}
If $u_{s_k}\in H^{s_k}_g(\Omega;B_L)$, with $s_k\to1^-$ as $k\to\infty$, is such that 
\begin{equation}
\label{eqi:s-jacs_liminf}
\jac^{s_{k}}(u_{s_{k}})\xflute{\overline\Omega}\pi\mu
\quad
\text{as}\
k\to\infty
\end{equation}
for some $\mu\in\atoms(\overline\Omega)$, then $\mu(\overline\Omega)=d$ and
\begin{equation*}
\liminf_{k\to\infty} 
\GL^{s_k}_\lambda(u_{s_k};\Omega)
\ge 
\pi|\mu|(\overline\Omega).
\end{equation*}

\item\label{item:main_limsup}
\emph{($\Gamma$-$\limsup$ inequality)}
If $\mu\in\atoms(\overline\Omega)$ is such that $\mu(\overline\Omega)=d$, then there exists a sequence $u_{s_k}\in H^{s_k}_g(\Omega;B_L)$, with $s_k\to1^-$ as $k\to\infty$, such that~\eqref{eqi:s-jacs_liminf} holds and
\begin{equation*}
\limsup_{k\to\infty} 
\GL^{s_k}_\lambda(u_{s_k};\Omega)
\le 
\pi|\mu|(\overline\Omega).\end{equation*}

\end{enumerate}
\end{theorem}

\subsection{Strategy of proof}
\label{subseci:proof}

The proof of \cref{resi:main} is split in two parts.

On the one side, claims~\ref{item:main_comp} and~\ref{item:main_liminf} of \cref{resi:main} follow by comparison with the local setting (recall \cref{resi:GL}) via a rigorous formulation of the approximation argument sketched in~\eqref{eqi:approx}.
The overall idea is to show that, if the fractional Ginzburg--Landau $s_k$-energy of $u_{s_k}$ is uniformly bounded, then  the integer $\e_k$-Ginzburg--Landau energy of $v_{\e_k}=I_{1-s_k}u_{s_k}$, with $\e_k=\sqrt{1-s_k}$, is also uniformly bounded and, actually, we have that
\begin{equation*}
\GL^{s_k}_\lambda(u_{s_k};\Omega)
\ge 
c_k
\GL_{\e_k,\Lambda}(v_{\e_k};\Omega)
\end{equation*}
for some $c_k>0$ such that $c_k\to\frac12$ as $k\to\infty$ and $\Lambda>0$ which does not depend on~$k$.
From \cref{resi:GL} we hence infer compactness and the $\Gamma$-$\liminf$ inequality for $\GL_{\lambda}^{s_k}(u_{s_k};\Omega)$ thanks to the fact that, due to the definition in~\eqref{eqi:def_s-jac} and~\eqref{eqi:riesz}, we have 
\begin{equation}
\label{eqi:jj}
\jac^{s_k}(u_{s_k})
=
c_k'\jac(v_{\e_k})
\end{equation} 
for some $c_k'>0$ such that $c_k'\to1$ as $k\to\infty$.
The fact that $\mu(\overline\Omega)=d$ follows from an integration-by-parts argument, roughly exploiting the fact that, recalling~\eqref{eqi:jj}, 
\begin{equation*}
\jac^{s_k}(u_{s_k})
=
c_k'
\jac(v_{s_k})
=
c_k'
\cu j(v_{s_k})
\sim
\cu j(g)
=
\jac(v_{s_k})
\end{equation*}
in a neighborhood $U$ of $\partial\Omega$, where $j(v_{s_k})=v_{s_k}\times Dv_{s_k}$, since 
\begin{equation*}
\|D^{s_k}u_{s_k}\|_{L^2(\R^2)}
\,
\|v_{s_k}-g\|_{L^2(U\setminus\Omega)}
\lesssim
\GL_\lambda^{s_k}(u_{s_k};\Omega)\,\sqrt{1-s_k}|\log(1-s_k)|
\to0^+,
\quad
\text{as}\
k\to\infty,
\end{equation*}
in view of the definition in~\eqref{eqi:s-GL} and of the approximation in~\eqref{eqi:approx} and~\eqref{eqi:approx2}. 

On the other side, concerning claim~\ref{item:main_limsup} in \cref{resi:main}, the recovery sequence $(u_{s_k})_{k\in\N}$ is built according to the classical cut-off approach around a vortex singularity.
Precisely, we first consider the case in which $\mu$ is given by
\begin{equation}
\label{eqi:special}
\mu=\sum_{i=1}^Nd_i\delta_{x_i},
\end{equation}
with $N\in\N$, $x_i\in\Omega$ and $d_i\in\set*{-1,1}$ such that 
\begin{equation}
\label{eqi:sum_deg}
\sum_{i=1}^Nd_i=\deg(g|_{\partial\Omega},\partial\Omega),
\end{equation}
and construct a recovery sequence of the form
\begin{equation}
\label{eqi:recovery}
u_{s_k}
=
\begin{cases}
\eta_k(\,\cdot\,-x_i)\,\dfrac{(\,\cdot\,-x_i)^{d_i}}{|\,\cdot\,-x_i|} & \text{around}\ x_i,\ \text{for each}\ i=1,\dots,N,
\\[1ex]
g 
& \text{outside}\ \Omega,
\\[1ex]
\widehat u
& \text{everywhere else},
\end{cases}
\end{equation}
where $\eta_k$ is a suitable cut-off near the origin depending on $s_k$, $x^d=\svec{x_1}{(-1)^d x_2}$ for $x\in\R^2$ and $d\in\set*{-1,1}$, and $\widehat u$ is an $H^1$ function with values in $\sphere$ which does not depend on $s_k$, but only on the position of the $x_i$'s in $\Omega$, interpolating the vortices $\frac{(\,\cdot\,-x_i)^{d_i}}{|\,\cdot\,-x_i|}$, $i=1,\dots,N$, and the datum $g$ (the existence of such a function $\widehat u$ is guaranteed by the condition~\eqref{eqi:sum_deg}).
Once the special case~\eqref{eqi:special} is achieved, the general situation can be then established via a routine diagonal approximation argument.   

\subsection{Comments}

We conclude this introduction with some comments on our main result.

An important difference between \cref{resi:GL} and \cref{resi:main} appears in the domains of definition of the corresponding energies. 
In our \cref{resi:main}, we only consider $H^s$ functions attaining the datum $g\in\bdata$ outside $\Omega$ which are bounded by some constant $L\ge\|g\|_{L^\infty}$.
We do not know if this additional constraint---which we need in the error estimate in~\eqref{eqi:approx2}---is merely technical and can be removed.
We nevertheless notice that such constraint does not seem to be essential in our approach.
Indeed, the fractional energies in~\eqref{eqi:s-GL} are decreased by truncation (see \cref{res:trunc} below) and the $L^\infty$ norm of the recovery sequence constructed in point~\ref{item:main_limsup} of \cref{resi:main} is bounded by~$\|g\|_{L^\infty}$. 

Another relevant point concerning the definition in~\eqref{eqi:sobs} is whether \cref{resi:main} may be achieved with respect to the smaller space
\begin{equation}
\label{eqi:sphereH}
H^s_g(\Omega;\sphere)
=
\set*{u\in H^s(\R^2;\R^2) : u=g\ \text{on}\ \R^2\setminus\Omega\ \text{and}\ |u|=1\ \text{on}\ \Omega}
\end{equation}
for $s\in(0,1)$, simplifying the non-local energies in~\eqref{eqi:s-GL} to
\begin{equation}
\label{eqi:s-nowell}
\GL^s(u;\Omega)
=
\frac{1}{|\log(1-s)|}
\int_{\R^2}|D^s u|\di x.
\end{equation}
While the compactness and the $\Gamma$-$\liminf$ follow as in the proof of \cref{resi:main} up to minor changes, the $\Gamma$-$\limsup$ remains unclear.
As done in~\cite{S23}, one may directly take $\frac{x^d}{|x|}$, with $x\ne0$ and $d\in\set*{1,-1}$, as a building block for the recovery sequence; that is, like in~\eqref{eqi:recovery}, but avoiding the truncation near the singularity.
However, the recovery sequence $\widetilde u_s\in H^1_g(\Omega;\sphere)$ constructed in this way satisfies
\begin{equation}
\label{eqi:s-fail}
\GL^s(\widetilde u_s;\Omega)
\lesssim
\frac{1}{(1-s)|\log(1-s)|}
\quad
\text{for}\
s\in(0,1).
\end{equation}
 
We do not know if the bound in~\eqref{eqi:s-fail} is optimal but, at least, it seems to agree with the one given by~\cite{S23}*{Th.~2.6} when (formally) rephrased to our setting.
Indeed, given a kernel $\rho\ge0$, the energy considered in~\cite{S23} is defined as 
\begin{equation}
\label{eqi:solci}
F_\e(u)
=
\frac{1}{\e^4|\log\e|}
\int_{\Omega}
\int_{\Omega}
\rho\Big(\frac{|x-y|}{\e}\Big)
|u(x)-u(y)|^2\di x\di y.
\end{equation}
Choosing $\rho(r)=r^{-2-2s}$ for $r>0$ (note that this kernel is not admissible for~\cite{S23}*{Th.~2.6}) and $\e=\sqrt{1-s}$, the functional in~\eqref{eqi:solci} becomes asymptotically equivalent as $s\to1^-$ to 
\begin{equation}
\label{eqi:s-solci}
u\mapsto
\frac{1}{|\log(1-s)|}
\int_{\Omega}
\int_{\Omega}
\frac{|u(x)-u(y)|^2}{|x-y|^{2+2s}}
\di x\di y.
\end{equation}
In virtue of~\eqref{eqi:s-energy} and~\eqref{eqi:s-equal}, and ignoring the reminder terms coming from the interactions between $\Omega$ and $\R^2\setminus\Omega$, the energy in~\eqref{eqi:s-solci} differs from the one in~\eqref{eqi:s-nowell} exactly by a factor $(1-s)$, which is precisely the additional factor appearing in the bound~\eqref{eqi:s-fail}.

The above considerations suggest that a more refined analysis of the energy~\eqref{eqi:s-fail} in the space in~\eqref{eqi:sphereH} is required in order to properly understand the $\Gamma$-limit
as $s\to 1^-$, and in particular whether vortex singularities are ruled out for our functionals with the
strict constraint $|u|=1$.

We finally note that functionals in~\eqref{eqi:s-GL} can be equivalently defined for $\Omega\subset\R^n$ with $n\ge 2$. 
In this case, \cref{resi:main} must be stated in terms of $(n-2)$-dimensional integral currents. 
We do not address the details here, but we refer to~\cite{ABO05} for proper statements and for the extension of the construction of recovery sequences for arbitrary $n\ge2$.

\subsection{Organization of the paper}
The main notation and definitions, plus some basic results concerning the Jacobian, the degree, the fractional gradient and the Riesz potential, are given in \cref{sec:preliminaries}.
The proof of \cref{resi:main} is then detailed across \cref{sec:proof}.

\section{Preliminaries}
\label{sec:preliminaries}

In this section, we recall the main notation and collect some preliminary results.

\subsection{General notation}
\label{subsec:notation}

Throughout the paper, the ambient space is~$\R^2$.
The norm of a point $x=\svec{x_1}{x_2}\in\R^2$ is given by $|x|^2=x_1^2+x_2^2$, while the norm of a matrix $A=\smat{A_{11}}{A_{12}}{A_{21}}{A_{22}}\in\R^{2\times2}$
is given by $|A|^2=A_{11}^2+A_{12}^2+A_{21}^2+A_{22}^2$.
We let $\sphere=\set*{x\in\R^2:|x|=1}$ be the standard unitary circle
and we let $\Omega\subset\R^2$ denote an open set, possibly satisfying further topological and/or regularity assumptions.

\subsection{Lipschitz functions, measures and flat norm}
\label{subsec:flat}

We let $\lip(\Omega)$ be the space of real-valued \emph{Lipschitz functions} on $\Omega$, endowed with the norm 
\begin{equation*}
\|\varphi\|_{\lip} 
=
\|\varphi\|_{L^\infty(\Omega)}
+
\|\nabla\varphi\|_{L^\infty(\Omega)},
\quad
\text{for}\ 
\varphi\in\lip(\Omega),
\end{equation*}
where $\nabla\varphi=\svec{\partial_{x_1}\varphi}{\partial_{x_2}\varphi}$ as customary,
and we let $\lip_c(\Omega)$ be its subspace of functions with compact support in $\Omega$.

We let $\mathcal M_b(\overline\Omega)$ be the space of \emph{Radon measures} on $\overline\Omega$ with finite total variation.
We let  
\begin{equation*}
\|\mu\|_{\flute(\Omega)}
=
\sup
\set*{
\int_\Omega\varphi\di\mu
:
\varphi\in\lip_c(\Omega),\ \|\varphi\|_{\lip}\le1
},
\end{equation*}
and, analogously,
\begin{equation*}
\|\mu\|_{\flute(\overline\Omega)}
=
\sup
\set*{
\int_\Omega\varphi\di\mu
:
\varphi\in\lip(\Omega),\ \|\varphi\|_{\lip}\le1
},
\end{equation*}
be the \emph{flat norm} of $\mu\in\mathcal M_b(\overline\Omega)$ on $\Omega$ and $\overline\Omega$, respectively.
Consequently, given $\mu_k,\mu\in\mathcal M_b(\overline\Omega)$, $k\in\N$, we write  
$\mu_k\xflute{\Omega}\mu$ (respectively,  
$\mu_k\xflute{\overline\Omega}\mu$)
if $\|\mu_n-\mu\|_{\flute(\Omega)}\to0$ (respectively, $\|\mu_k-\mu\|_{\flute(\overline\Omega)}\to0$) as $k\to\infty$.

\subsection{Sobolev functions, Jacobian and degree}
\label{subsec:sobolev}

We let  
\begin{equation*}
H^1(\Omega;\R^2)
=
\set*{v\in L^2(\Omega;\R^2) : Dv\in L^2(\Omega;\R^{2\times 2})}
\end{equation*}
be the standard \emph{Sobolev space} on $\Omega$, endowed with the norm
\begin{equation*}
\|v\|_{H^1}
=
\|v\|_{L^2(\Omega)}+\|Dv\|_{L^2(\Omega)},
\quad
\text{for}\ 
v\in H^1(\Omega;\R^2),
\end{equation*}
see~\cite{L17book}*{Ch.~11} for a detailed introduction.
Here we note that
\begin{equation*}
Dv=\mat{\partial_{x_1}v_1}{\partial_{x_2}v_1}{\partial_{x_1}v_2}{\partial_{x_2}v_2}
\end{equation*}
and that $\|Dv\|_{L^2(\Omega)}^2
=
\||Dv|\|_{L^2(\Omega)}^2$, where $|Dv|$ is defined as in \cref{subsec:notation}.

The \emph{Jacobian} of $v\in H^1(\Omega;\R^2)$ is defined as 
\begin{equation}
\label{eq:jac}
\jac(v)=\det(Dv)
=
\partial_{x_1}v_1
\partial_{x_2}v_2
-
\partial_{x_2}v_1
\partial_{x_1}v_2
\in 
L^1(\Omega).
\end{equation}
In particular, we tacitly identify $\jac(v)$ with the Radon measure $\jac(v)\mathscr L^2\in\mathcal M_b(\overline\Omega)$.
For a (historical) presentation of the Jacobian of Sobolev functions, we refer to the survey~\cite{BMM24}.

We now recall some elementary properties of the Jacobian that will be useful in the sequel.
The first one concerns $H^1$ functions with values in~$\sphere$ (e.g., see~\cite{BMM24}*{Sec.~3}). 

\begin{lemma}
\label{res:null_jac}
If $v\in H^1(\Omega;\R^2)$ is such that $|v|=1$ on $\Omega$, then $\jac(v)=0$ on $\Omega$.
\end{lemma}

\begin{proof}
Differentiating $|v|^2=1$, we get that $v\in\ker(Du)$, proving that $\det(Dv)=0$.
\end{proof}

In order to state the second result on the Jacobian, we need the following notation. 

\begin{definition}
\label{def:j}
Given $v,w\in H^1(\Omega;\R^2)$, we define the \emph{current}  
\begin{equation*}
j(v,w)=v\times Dw
=
v_1Dw_2-v_2Dw_1
\in L^1(\Omega;\R^2).
\end{equation*}
In particular, we let $j(v)=j(v,v)\in L^1(\Omega;\R^2)$.
\end{definition}

The next result relates the Jacobian of $H^1$ functions to the current $j$ in \cref{def:j} and also shows that the $\cu$ of $j$ is zero. 

\begin{lemma}
\label{res:cucu}
If $v,w\in H^1(\Omega;\R^2)$, then 
\begin{equation}
\label{eq:cucu_jac}
\int_\Omega\jac(v)\,\varphi\di x
=
\frac12
\int_\Omega j(v)\times\nabla\varphi\di x
\end{equation}
and
\begin{equation}
\label{eq:cucu_sym}
\int_\Omega (j(v,w)-j(w,v))\times \nabla\varphi\di x
=0
\end{equation}
whenever $\varphi\in\lip_c(\Omega)$.
\end{lemma}

\begin{proof}
If $v,w\in C^2(\Omega;\R^2)$, then we can write 
\begin{equation*}
\jac(v)=\frac12\cu j(v)
\quad
\text{and}
\quad
j(v,w)-j(w,v)
=
\nabla(v\times w),
\end{equation*}
from which we get~\eqref{eq:cucu_jac} and~\eqref{eq:cucu_sym}.
In the general case, the validity of~\eqref{eq:cucu_jac} and~\eqref{eq:cucu_sym} follows by approximating $v,w\in H^1(\Omega;\R^2)$ with smooth functions.
\end{proof}

\cref{res:cucu} can be exploited to estimate the flat distance between the Jacobians of two $H^1$ functions, as follows (note that this is a particular instance of~\cite{BN11}*{Th.~1(i)}).

\begin{proposition}
\label{res:flat_dist_jac}
If $v,w\in H^1(\Omega;\R^2)$, then
\begin{equation*}
\|\jac(v) - \jac(w)\|_{\flute(\Omega)}
\le 
2\,\|v - w\|_{L^2(\Omega)}\,(\|Dv\|_{L^2(\Omega)} + \|Dw\|_{L^2(\Omega)}).
\end{equation*}
\end{proposition}

\begin{proof}
Given $v,w\in H^1(\Omega;\R^2)$, by \cref{def:j} we can estimate
\begin{equation*}
|j(v)-j(w)|
\le 
|j(v,v-w)|
+
|j(w,v-w)|
\le 
4\,|v-w|\,(|Dv|+|Dw|).
\end{equation*}
From~\eqref{eq:cucu_jac}, by H\"older's inequality we hence get that 
\begin{equation*}
\left|\int_\Omega(\jac(v)-\jac(w))\,\varphi\di x
\right|
\le 
2\,\|v - w\|_{L^2(\Omega)}\,(\|Dv\|_{L^2(\Omega)} + \|Dw\|_{L^2(\Omega)})\|\nabla\varphi\|_{L^\infty(\Omega)}
\end{equation*}
whenever $\varphi\in\lip_c(\Omega)$, yielding the conclusion.
\end{proof}

Last, but not least, we recall the following result linking the Jacobian of Sobolev functions to their (\emph{topological}) \emph{degree}.
For a more detailed account, we refer to the survey~\cite{B97}.

\begin{lemma}
\label{res:deg}
Let $\Omega\subset\R^2$ be a non-empty, connected, simply connected, bounded open set with Lipschitz boundary.
If $u\in H^1(\Omega;\R^2)$ is such that $|u|=1$ on $\partial\Omega$, then the (topological) degree $\deg(u|_{\partial\Omega},\partial\Omega)\in\Z$ of the trace $u|_{\partial\Omega}$ of $u$ on $\partial\Omega$ satisfies
\begin{equation}
\label{eq:deg}
\deg(u|_{\partial\Omega},\partial\Omega)
=
\frac1{\pi}
\int_\Omega \jac(u) \di x.
\end{equation} 
\end{lemma}

\subsection{Ginzburg--Landau energies}

As briefly recalled in \cref{subseci:classical}, given $\e,\lambda>0$, the \emph{Ginzburg--Landau functionals} on a non-empty open set $\Omega$,
\begin{equation*}
\GL_{\e,\lambda}(\,\cdot\,;\Omega)\colon H^1(\Omega;\R^2)\to\R,
\end{equation*}
are defined as 
\begin{equation}
\label{eq:GLe}
\GL_{\e,\lambda}(u;\Omega)
= 
\frac1{|\log\e|}
\int_{\Omega}|D u|^2\di x
+
\frac{\lambda}{\e^2|\log \e|}\int_{\Omega}(|u|^2-1)^2\di x
\end{equation}
for $u\in H^1(\Omega;\R^2)$.
The following result rephrases \cref{resi:GL} in the case no boundary condition is imposed on $\partial\Omega$ (again, refer to~\cites{J99,JS02,S00,S98,BBH94,AP14} for the proof, and see~\cite{ABO05} for the higher-dimensional setting).
Here and in below, as in~\eqref{eqi:atoms}, we let
\begin{equation*}
\atoms(A)
=
\set*{\mu=\sum_{i=1}^N d_i\,\delta_{x_i} : d_i\in\Z\ \text{and}\ x_i\in A\ \text{for}\ i=1,\dots, N,\ \text{with}\ N\in\N}
\end{equation*}
be the collections of \emph{atomic} measures on a non-empty set $A\subset\R^2$.

\begin{theorem}
\label{res:GL}
Let $\Omega\subset\R^2$ be a non-empty, connected, simply connected, bounded open set with Lipschitz boundary and $\lambda>0$.

\begin{enumerate}[label=(\roman*),ref=(\roman*),itemsep=1ex,leftmargin=4ex]

\item
\label{item:GL_compactness}
\emph{(Compactness)}
If $(v_{\e_k})_{k\in\N}\subset H^1(\Omega;\R^2)$, with $\e_k\to0^+$ as $k\to\infty$, is such that 
\begin{equation*}
\sup_{k\in\N}\GL_{\e_k,\lambda}(v_{\e_k};\Omega)<\infty,
\end{equation*}
then there exists a subsequence $(v_{\e_{k_j}})_{j\in\N}$ and $\mu\in\atoms(\overline\Omega)$ such that
\begin{equation*}
\jac(v_{\e_{k_j}})
\xflute{\overline\Omega}\mu
\quad
\text{as}\ j\to\infty.
\end{equation*}

\item
\label{item:GL_liminf}
\emph{($\Gamma$-$\liminf$ inequality)}
If $(v_{\e_k})_{k\in\N}\subset H^1_g(\Omega;\R^2)$, with $\e_k\to0^+$ as $k\to\infty$, is such that
\begin{equation}
\label{eq:jacs_liminf}
\jac(v_{\e_k})
\xflute{\overline\Omega}\mu
\quad
\text{as}\ k\to\infty
\end{equation}
for some $\mu\in\atoms(\overline\Omega)$, then 
\begin{equation*}
\liminf_{k\to\infty}\GL_{\e_k,\lambda}(v_{\e_k};\Omega)
\ge 
2\pi|\mu|(\overline\Omega).
\end{equation*}

\item
\label{item:GL_limsup}
\emph{($\Gamma$-$\limsup$ inequality)}
If $\mu\in\atoms(\overline\Omega)$, then there exists a sequence $(v_{\e_k})_{k\in\N}\subset H^1_g(\Omega;\R^2)$, with $\e_k\to0^+$ as $k\to\infty$, such that~\eqref{eq:jacs_liminf} holds and 
\begin{equation*}
\limsup_{k\to\infty}\GL_{\e_k,\lambda}(v_{\e_k};\Omega)
\le 
2\pi|\mu|(\overline\Omega).
\end{equation*}  

\end{enumerate}
\end{theorem}

\subsection{Fractional Sobolev functions}

We let 
\begin{equation*}
H^s(\R^2;\R^2)
=
\set*{u\in L^2(\Omega;\R^2) : [u]_{s,\Omega}<\infty}
\end{equation*}
be the \emph{fractional Sobolev space} of order $s\in(0,1)$, endowed with the norm
\begin{equation*}
\|u\|_{H^s}
=
\|u\|_{L^2}
+
[u]_{s,\Omega},
\quad
\text{for}\ u\in H^s(\Omega;\R^2),
\end{equation*}
where 
\begin{equation}
\label{eq:seminorm_Hs}
[u]_{s,\Omega}^2
=
\int_{\Omega}\int_{\Omega}\frac{|u(y)-u(x)|^2}{|y-x|^{2+2s}}\di x\di y,
\end{equation}
see~\cite{L23book}*{Ch.~6} for a detailed introduction.
If $\Omega=\R^2$, then we simply write $[u]_s=[u]_{s,\R^2}$.

For future convenience, we recall the following result, which corresponds to~\cite{BBD23}*{Lem.~3}.
Here and below, for $\tau>0$, we let 
\begin{equation}
\label{eq:bermuda}
\bigtriangleup_\tau
=
\set*{(x,y)\in\R^2 : |x-y|\le\tau}
\subset\R^2.
\end{equation}

\begin{lemma}
\label{res:fartau}
Let $\Omega\subset\R^2$ be a measurable set.
If $u\in L^2(\Omega)$, then 
\begin{equation*}
\iint_{(\Omega\times\Omega)\setminus\bigtriangleup_\tau}
\frac{|u(x)-u(y)|^2}{|x-y|^{2+2s}}\di x\di y
\le
\frac{4\pi\|u\|_{L^2(\Omega)}^2}{s\tau^{2s}}
\end{equation*}
for $\tau>0$ and $s\in(0,1)$.
\end{lemma}

\begin{proof}
By the definition in~\eqref{eq:bermuda}, we can estimate
\begin{equation*}
\iint_{(\Omega\times\Omega)\setminus\bigtriangleup_\tau}
\frac{|u(x)-u(y)|^2}{|x-y|^{2+2s}}\di x\di y
\le 
4
\int_{\Omega}|u(x)|^2
\int_{\R^2\setminus B_\tau}
\frac{\di h}{|h|^{2+2s}}\di x\di y
=
\frac{4\pi\|u\|_{L^2(\Omega)}^2}{s\tau^{2s}}
\end{equation*}
whenever $\tau>0$ and $s\in(0,1)$.
\end{proof}

\subsection{Fractional gradient}

By combining~\cite{SS15}*{Rem.~2.3} with~\cite{BCCS22}*{Cor.~1} (also see the discussion in~\cite{CS19}*{Sec.~3.9}), we can equivalently define 
\begin{equation*}
H^s(\R^2;\R^2)
=
\set*{u\in L^2(\R^2;\R^2) : D^su\in L^2(\R^2;\R^{2\times 2})}
\end{equation*}
for $s\in(0,1)$, where, as briefly recalled in~\eqref{eqi:s-grad}, 
\begin{equation}
\label{eq:s-grad}
D^s u(x)
=
(1-s)\,\dfrac{2^{s-1}}{\pi}\,\tfrac{\Gamma\left(\frac{3+s}{2}\right)}{\Gamma\left(\frac{3-s}{2}\right)}
\int_{\R^2}\frac{(u(y)-u(x))\otimes(y-x)}{|y-x|^{3+s}}\di y,
\quad
x\in\R^2,
\end{equation}
is the \emph{fractional} (\emph{Riesz}) \emph{$s$-gradient} of ~$u$.
Note that the $s$-gradient in~\eqref{eq:s-grad} is well defined for sufficiently regular functions ($u\in\lip_c(\R^2;\R^2)$ would suffice, see the discussion in~\cite{CS19}*{Sec.~2.2} for instance), while, for $u\in H^s(\R^2;\R^2)$, the operator in~\eqref{eq:s-grad} is defined in the distributional sense via (fractional) integration-by-parts, see~\cite{CS19}*{Def.~3.19}.  

Here we just recall the following result, which relates the $L^2$ norm of the $s$-gradient in~\eqref{eq:s-grad} with the fractional seminorm in~\eqref{eq:seminorm_Hs}.
Here and below, as in \cref{subsec:sobolev}, we set $\||D^su\|_{L^2}^2=\||D^su|\|_{L^2}^2$, where $|D^su| $ is as defined in \cref{subsec:notation}.

\begin{proposition}
\label{res:equiv-frac}
If $u\in H^s(\mathbb R^2;\mathbb R^2)$, then 
\begin{equation}
\label{eq:relation_s-grad_seminorm}
\|D^s u\|_{L^2}^2
=
(1-s)\,c_s\,[u]^2_s, 
\end{equation}
where $c_s=\dfrac{4^s}{2\pi}\,\dfrac{s\,\Gamma(1+s)}{\Gamma(2-s)}>0$ satisfies $\lim\limits_{s\to 1^-}c_s=\dfrac{2}{\pi} $.
\end{proposition}

\begin{proof}
By density, we may assume that $u\in C^\infty_c(\R^2;\R^2)$ without loss of generality.
In this case, formula~\eqref{eq:relation_s-grad_seminorm} follows either by applying Fourier's transform, or by observing that 
\begin{equation*}
(1-s)\,c_s\,[u]_s^2
=
\int_{\R^2}u\,(-\Delta)^su\di x
=
-\int_{\R^2}u\,\diverg^s(D^su)\di x
=
\int_{\R^2}|D^s u|^2\di x,
\end{equation*}
where $(-\Delta)^s$ is the \emph{fractional Laplacian} of order $s$ and $\diverg^s$ is the \emph{fractional divergence} of order $s$ (the dual operator of~\eqref{eq:s-grad}), see~\cite{CS19}*{Lem.~2.5 and Sec.~3.10}.
\end{proof}

From \cref{res:equiv-frac}, we deduce that the fractional Ginzburg--Landau energies in~\eqref{eqi:s-GL} are decreased by truncation.
In order to precisely state this result, we need to introduce some notation. 
Given $L>0$, we let 
\begin{equation*}
T_L(t)=\max\set*{-L,\min\set*{t,L}},
\quad
\text{for}\ t\in\R,
\end{equation*}
and we set
\begin{equation*}
T_L(x)=\vect{T_L(x_1)}{T_L(x_2)},
\quad
\text{for}\
x=\vect{x_1}{x_2}\in\R^2.
\end{equation*} 
Note that $|T_L(x)|\le|x|$ for all $x\in\R^2$ by definition.

\begin{corollary}
\label{res:trunc}
If $u\in H^s(\R^2;\R^2)$, then $T_L(u)\in H^s(\R^2;\R^2)$ and 
\begin{equation}
\label{eq:trunc}
\GL^s_{\lambda}(T_L(u))
\le
\GL^s_{\lambda}(u)
\end{equation}
for all $s\in(0,1)$ and $L\ge1$.
\end{corollary}

\begin{proof}
Since $T_L\colon\R\to\R$ is $1$-Lipschitz, from \cref{res:equiv-frac} we get that
\begin{equation*}
\|D^s(T_L(u))\|^2_{L^2}
=
(1-s)\,c_s\,[T_L(u)]_s^2
\le 
(1-s)\,c_s\,[u]_s^2
=
\|D^s u\|^2_{L^2}
\end{equation*}
whenever $L>0$. 
Moreover, for $L\ge1$, we also have that $(|T_L(u)|^2-1)^2\le(|u|^2-1)^2$, from which the inequality in~\eqref{eq:trunc} follows.
\end{proof}

\subsection{Riesz potential}

As observed in~\cite{SS15}*{Th.~1.2} (also see~\cite{CS19}*{Sec.~2.3}), the $s$-gradient in~\eqref{eq:s-grad} can be rewritten as $D^s=D I_{1-s}$, where
\begin{equation}
\label{eq:riesz}
I_{1-s}u(x)
=
\frac{1}{\gamma_s}
\int_{\R^2}\frac{u(y)}{|y-x|^{1+s}}\di y,
\quad
x\in\R^2,
\end{equation}
is the \emph{Riesz potential} of order $1-s$ in $\R^2$, with
\begin{equation}
\label{eq:gamma_s}
\gamma_s
=
2^{1-s}\pi
\,
\frac{1+s}{1-s}
\,
\tfrac{\Gamma\left(\frac{3-s}{2}\right)}{\Gamma\left(\frac{3+s}{2}\right)}.
\end{equation}
Precisely, considering the differential operators in the distributional sense, we have that
\begin{equation}
\label{eq:nabla_of_riesz}
D^s u=D I_{1-s}u\
\text{in}\ 
L^2(\R^2;\R^{2\times2})
\end{equation}
for $u\in H^s(\R^2;\R^2)$.
Therefore, if $u\in H^s(\R^2;\R^2)$, then 
\begin{equation}
\label{eq:riesz_sob}
v=I_{1-s}u\in H^1(\R^2;\R^2)
\end{equation}
and so, recalling~\eqref{eq:jac}, we may generalize the notion of Jacobian to fractional Sobolev functions as follows.

\begin{definition}[Fractional Jacobian]
\label{def:frac_jac}
Given $s\in(0,1)$, the \emph{fractional} (\emph{Riesz}) \emph{$s$-Jacobian} of $u\in H^s(\R^2;\R^2)$ is defined as 
\begin{equation*}
\begin{split}
\jac^s(u)
=
\jac(I_{1-s}u)
\in L^1(\R^2).
\end{split}
\end{equation*}
\end{definition}

In particular, we tacitly identify $\jac^s(u)$ with the Radon measure $\jac^s(u)\mathscr L^2\in\mathcal M_b(\R^2)$ whenever $s\in(0,1)$ and $u\in H^s(\R^2;\R^2)$.

\begin{remark}
Although not needed in the present work, \cref{def:j,res:cucu,res:flat_dist_jac} can be naturally analogously reformulated in the fractional setting.
\end{remark}

\section{Proof of \texorpdfstring{\cref{resi:main}}{the main result}}
\label{sec:proof}

The rest of the paper is dedicated to the proof of \cref{resi:main}, which is split across \cref{subsec:compactness,subsec:liminf,subsec:limsup}.
From now on, with the notation introduced in~\eqref{eqi:bdata}, we fix a boundary datum $g\in\bdata$ and a bounded open neighborhood $U\subset\R^2$ of $\partial\Omega$ such that $|g|=1$ on~$U$.
We hence let $L\in [\|g\|_{L^\infty(\R^2)},\infty )$, we let $A=U\cup\Omega$ and $R>0$ be such that $\supp \bar u\subset A_R$, where
\begin{equation*}
A_R=\bigcup_{x\in A}B_R(x),
\end{equation*}
and we work in the fractional Sobolev space $H^s_g(\Omega;B_L)$ defined in~\eqref{eqi:sobs}.

\subsection{A truncated Riesz potential}
\label{subsec:truncated_Riesz}

We begin with the following definition, which provides a truncated version of the Riesz potential defined in~\eqref{eq:riesz}.

\begin{definition}[Truncated Riesz potential]
\label{def:truncated_Riesz}
For $s\in(0,1)$, the \emph{truncated Riesz potential} $I^R_{1-s} u\colon\R^2\to\R^2$ of $u\in H^s_g(\Omega;B_L)$ is defined as
\begin{equation}
\label{eq:truncated_Riesz}
I^R_{1-s} u(x)
=
\frac{1-s}{2\pi R^{1-s}}
\int_{B_R(x)}\frac{u(y)}{|x-y|^{1+s}}\di y,
\quad
\text{for}\ x\in\R^2.
\end{equation}
\end{definition}

Let  $u\in H^s_g(\Omega;B_L)$ be fixed, and set $v=I^R_{1-s}u$ as in \cref{def:truncated_Riesz} above for brevity.
Recalling~\eqref{eq:riesz}, and owing to the fact that $\supp g\subset A$, we can write
\begin{equation}
\label{confront}
v(x)
=
\frac{(1-s)\gamma_s}{2\pi R^{1-s}}
\,
I_{1-s}u(x),
\quad
\text{for}\
x\in A,
\end{equation}
with $\gamma_s>0$ as in~\eqref{eq:gamma_s} such that 
\begin{equation}
\label{eq:lim_gamma_R}
\lim_{s\to1^-}
\frac{(1-s)\gamma_s}{2\pi R^{1-s}}
=
1.
\end{equation}
Moreover, by combining~\eqref{eq:nabla_of_riesz} and~\eqref{eq:riesz_sob} with~\eqref{confront}, we get that $v\in H^1(A;B_M)$, with
\begin{equation}
\label{gradok}
D v
=
\frac{(1-s)\gamma_s}{2\pi R^{1-s}}\,D^su
\quad
\text{in}\ L^2(A;\R^{2\times2})
\end{equation}
and 
\begin{equation}
\label{eq:jacs}
\jac^s(u)
=
\left(\frac{2\pi R^{1-s}}{(1-s)\gamma_s}\right)^2
\jac(v)
\quad
\text{on}\ A.
\end{equation}

We now exploit \cref{def:truncated_Riesz} in two ways.
We first provide a quantitative formulation of the error in the approximation argument carried in~\eqref{eqi:approx} and~\eqref{eqi:approx2}.

\begin{lemma}[$L^2$ comparison]
\label{res:L2_comp}
If $u\in\sobs$ then $v=I_{1-s}^Ru$ as in \cref{def:truncated_Riesz} satisfies
\begin{equation}
\label{eq:L2_comp}
\int_{A}\left(\left|v\right|^2-|u|^2\right)^2\di x
\le
4L^2 \int_A \left|u-v\right|^2\di x
\le 
\frac{(1-s)^2 L^2}{\pi R^{1-2s}}
\,[u]_s^2
\end{equation}
for any $s\in(0,1)$.
\end{lemma}

\begin{proof}
Since
\begin{equation*}
u(x)
=
\frac{1-s}{2\pi R^{1-s}}
\int_{B_R(x)}\frac{u(x)}{|x-y|^{1+s}}\di y,
\quad
\text{for}\ x\in\R^2,
\end{equation*}
we can estimate
\begin{equation*}
|u+v|^2
=
\bigg|\frac{1-s}{2\pi R^{1-s}}\int_{B_R(x)}\frac{u(x)+u(y)}{|x-y|^{1+s}}\di y\bigg|^2
\le 
4\|u\|_{L^\infty(A_R)}^2
\le 
4L^2 
\end{equation*}
for any $x\in A$,
so that
\begin{equation*}
\begin{split}
\left(\left|v(x)\right|^2-|u(x)|^2\right)^2
&\le 
\left|v(x) +u(x)\right|^2
\left|v(x) -u(x)\right|^2
\le 
4L^2
\left|v(x) -u(x)\right|^2
\end{split}
\end{equation*}
for any $x\in A$, giving the first inequality in~\eqref{eq:L2_comp} by integrating on $A$.
For the second inequality in~\eqref{eq:L2_comp}, we just need to observe that, by Jensen's inequality,
\begin{equation*}
\begin{split}
\left|v(x)-u(x)\right|^2
=
\left|\frac{1-s}{2\pi R^{1-s}}\int_{B_R(x)}\frac{u(x)-u(y)}{|x-y|^{1+s}}\di y\right|^2
\le
\frac{(1-s)^2}{\pi R^{1-2s}} \int_{B_R(x)}\frac{|u(x)-u(y)|^2}{|x-y|^{2+2s}}\di y
\end{split}
\end{equation*}
for any $x\in A$, and the conclusion  follows by integrating on~$A$ again. 
\end{proof}

Additionally, we compare the Ginzburg--Landau energies in~\eqref{eq:GLe} with their fractional counterparts in~\eqref{eqi:s-GL}.
As mentioned in \cref{subseci:proof}, \cref{res:E_comp} below will play a crucial role in the proof of claims~\ref{item:main_comp} and~\ref{item:main_liminf} of \cref{resi:main}.

\begin{lemma}[Energy comparison]
\label{res:E_comp}
If $u\in\sobs$, then $v=I_{1-s}^Ru$ as in \cref{def:truncated_Riesz} satisfies
\begin{align}
\GL^s_{\lambda}(u;\Omega)
&\ge 
\frac{4\pi R^{1-2s}c_s}{(1-s)|\log(1-s)|}
\int_{A}|u-v|^2\di x,
\label{eq:E_comp_L2}
\\[1ex]
\GL^s_{\lambda}(u;\Omega)
&\ge 
\frac{\sqrt{4\pi R^{1-2s}c_s}}{\sqrt{1-s}\,|\log(1-s)|}
\,
\|D^s u\|_{L^2(\R^2)}\,\|v-u\|_{L^2(A)},
\label{eq:E_comp_prod}
\\[1ex]
\GL^s_{\lambda}(u;\Omega)
&\ge 
\frac{\eta}2 
\left(\frac{2\pi R^{1-s}}{(1-s)\gamma_s}\right)^2
\,
\GL_{\sqrt{1-s},\,\Lambda_{s,\lambda,\eta}}(v;A),
\label{eq:E_comp_GL}
\end{align}
for $s\in(0,1)$, $\lambda>0$ and $\eta\in(0,1)$, where  
\begin{equation}
\label{eq:E_comp_const}
\Lambda_{s,\lambda,\eta}
=
\frac{1-\eta}{2\eta}
\,
\left(\frac{(1-s)\gamma_s}{2\pi R^{1-s}}\right)^2
m_{s,\lambda},
\quad
m_{s,\lambda}
=
\min\set*{\frac{c_s\pi R^{1-2s}}{L^2},\lambda},
\end{equation}
and $c_s,\gamma_s>0$ are as in~\eqref{eq:relation_s-grad_seminorm} and~\eqref{eq:gamma_s}, respectively. 
\end{lemma}

\begin{proof}
By \cref{res:equiv-frac} and the first inequality in~\eqref{eq:L2_comp} in \cref{res:L2_comp}, we get that 
\begin{equation*}
\GL_\lambda^s(u;\Omega)
\ge 
\frac{(1-s)\,c_s}{|\log(1-s)|}
\,
[u]_s^2
\ge 
\frac{4\pi R^{1-2s}c_s}{(1-s)|\log(1-s)|}
\int_{A}|u-v|^2\di x,
\end{equation*}
proving~\eqref{eq:E_comp_L2}.
Combining~\eqref{eq:E_comp_L2} with the inequality
\begin{equation*}
\GL_\lambda^s(u;\Omega)
\ge 
\frac{1}{|\log(1-s)|}
\int_{\R^2}|D^s u|^2\di x,
\end{equation*}
we get~\eqref{eq:E_comp_prod}.
Owing to the fact that $u=g$ and $|g|=1$ on $A\setminus\Omega$, \cref{res:equiv-frac} and the second inequality in~\eqref{eq:L2_comp} in \cref{res:L2_comp}, we can estimate
\begin{equation}
\label{scendo}
\begin{split}
\GL_\lambda^s(u;\Omega)
&=
\frac{(1-s)\,c_s}{|\log(1-s)|}
\,
[u]_s^2
+
\frac{\lambda}{(1-s)|\log(1-s)|}
\int_{A}\left(|u|^2-1\right)^2
\di x
\\
&\ge
\frac{(1-s)\,c_s}{|\log(1-s)|}
\,
\frac{\pi R^{1-2s}}{(1-s)^2 L^2}
\int_{A}\left(|v|^2-|u|^2\right)^2\di x
\\
&\quad
+
\frac{\lambda}{(1-s)|\log(1-s)|}\int_{A}\left(|u|^2-1\right)^2\di x
\\
&\ge
\frac{m_{s,\lambda}}{2\,(1-s)|\log(1-s)|}
\int_{A}\left(|v|^2-1\right)^2\di x.
\end{split}
\end{equation}
Therefore, by~\eqref{scendo} and~\eqref{gradok}, we get that 
\begin{equation*}
\begin{split}
&\GL_\lambda^s(u)
=
\eta\GL_\lambda^s(u)+(1-\eta)\GL_\lambda^s(u)
\\
&\quad\ge
\frac{\eta}{|\log(1-s)|}\int_A|D^s u|^2\di x
+
\frac{(1-\eta)\,m_{s,\lambda}}{2\,(1-s)|\log(1-s)|}
\int_{A}\left(|v|^2-1\right)^2\di x
\\
&\quad\ge
\frac{\eta}2 
\left(\frac{2\pi R^{1-s}}{(1-s)\gamma_s}\right)^2
\frac{1}{|\log\sqrt{1-s}|}
\int_A
|D v|^2\di x
+
\frac{(1-\eta)\,m_{s,\lambda}}{4\,(1-s)|\log\sqrt{1-s}|}
\int_{A}\left(|v|^2-1\right)^2\di x,
\end{split}
\end{equation*} 
for any $\eta\in(0,1)$, from which~\eqref{eq:E_comp_GL} follows due to the definition in~\eqref{eq:GLe}.
\end{proof}

For future convenience, completing the definitions in~\eqref{eq:E_comp_const} in \cref{res:E_comp}, we set
\begin{equation}
\label{eq:lim_Lambda}
\Lambda_{1,\lambda,\eta}
=
\lim_{s\to1^-}
\Lambda_{s,\lambda,\eta}
=
\frac{1-\eta}{2\eta}\,m_{1,\lambda},
\end{equation}
for any $\lambda>0$ and $\eta\in(0,1)$, where
\begin{equation}
\label{eq:lim_c_R}
m_{1,\lambda}
=
\lim_{s\to1^-}
m_{s,\lambda}
=
\min\set*{\frac{1}{RL^2},\lambda}.
\end{equation}

\subsection{Proof of claim \texorpdfstring{\ref{item:main_comp}}{the compactness} of \texorpdfstring{\cref{resi:main}}{the main result}}
\label{subsec:compactness}

Let
$u_{s_k}\in H^{s_k}_g(\Omega;B_L)$, with $s_k\to1^-$ as $k\to\infty$, be such that 
\begin{equation}
\label{eq:sup_F}
\sup_{k\in\N} 
\GL_\lambda^{s_k}(u_{s_k};\Omega)<\infty.
\end{equation}
Following \cref{def:truncated_Riesz}, we define $v_{s_k}=I_{1-{s_k}}^Ru_{s_k}$ for $k\in\N$.
By combining~\eqref{eq:lim_gamma_R} and~\eqref{eq:lim_c_R} with~\eqref{eq:E_comp_GL} in \cref{res:E_comp}, from~\eqref{eq:sup_F} we infer that
\begin{equation*}
\sup_{k\in\N}\GL_{\sqrt{1-{s_k}},\,C}(v_{s_k};A)<\infty
\end{equation*}
for some $C>0$ which does not depend on~$k$. 
Thus, by \cref{res:GL}\ref{item:GL_compactness}, we find $\mu\in\atoms(A)$ such that, up to passing to a subsequence (which we do not relabel), 
\begin{equation}
\label{eq:conv_int_jac}
\jac(v_{s_k})\xflute{A}\pi\mu
\quad
\text{as}\
k\to\infty.
\end{equation}
Thanks to~\eqref{eq:jacs}, we can  estimate
\begin{equation*}
\begin{split}
\|\jac^{s_k}(u_{s_k})
-\pi\mu\|_{\flute(A)}
\le 
\left|1-\left(\frac{2\pi R^{1-{s_k}}}{(1-{s_k})\gamma_{s_k}}\right)^2\right|
\,
\|\jac(v_{s_k})\|_{\flute(A)}
+
\|\jac(v_{s_k})
-\pi\mu\|_{\flute(A)},
\end{split}
\end{equation*}
so that, by~\eqref{eq:lim_gamma_R} and~\eqref{eq:conv_int_jac}, we deduce that
\begin{equation}
\label{eq:jac_conv_large}
\jac^{s_k}(u_{s_k})\xflute{A}\pi\mu
\quad
\text{as}\
k\to\infty.
\end{equation}
Since $\Omega\Subset A$ by construction, from~\eqref{eq:jac_conv_large} in particular we get that
\begin{equation*}
\jac^{s_k}(u_{s_k})\xflute{\overline\Omega}\pi\mu
\quad
\text{as}\
k\to\infty.
\end{equation*}
We now prove that $\|\mu\|_{\flute(A\setminus\Omega)}=0$, so that $\supp \mu\subset\overline\Omega$. 
To this aim, we observe that
\begin{equation}
\label{eq:salire}
\pi
\|\mu\|_{\flute(A\setminus\Omega)}
\le 
\|\jac^{s_k}(u_{s_k})-\pi\mu\|_{\flute(A)}
+
\|\jac^{s_k}(u_{s_k})\|_{\flute(A\setminus\Omega)}.
\end{equation}
In view of~\eqref{eq:jac_conv_large}, we just need to deal with the second term in the right-hand side.
Owing to~\eqref{eq:jacs}, \cref{res:null_jac} (since $|g|=1$ on $A\setminus\Omega$), \cref{res:flat_dist_jac} and~\eqref{gradok}, we have that
\begin{equation}
\label{eq:mischia}
\begin{split}
\|\jac^{s_k}&(u_{s_k})\|_{\flute(A\setminus\Omega)}
=
\left(\tfrac{2\pi R^{1-{s_k}}}{(1-{s_k})\gamma_s}\right)^2
\|\jac(v_s)\|_{\flute(A\setminus\Omega)}
\\
&=
\left(\tfrac{2\pi R^{1-{s_k}}}{(1-{s_k})\gamma_{s_k}}\right)^2
\|\jac(v_{s_k})-\jac(g)\|_{\flute(A\setminus\Omega)}
\\
&\le 
C\left(\tfrac{2\pi R^{1-{s_k}}}{(1-{s_k})\gamma_{s_k}}\right)^2
\|v_{s_k}-g\|_{L^2(A\setminus\Omega)}
\,
\left(
\|Dv_{s_k}\|_{L^2(A\setminus\Omega)}
+
\|Dg\|_{L^2(A\setminus\Omega)}
\right)
\\
&\le 
C\left(\tfrac{2\pi R^{1-{s_k}}}{(1-{s_k})\gamma_{s_k}}\right)^2
\|v_{s_k}-g\|_{L^2(A\setminus\Omega)}
\,
\left(
\tfrac{(1-{s_k})\gamma_s}{2\pi R^{1-s}}
\,
\|D^{s_k} u_{s_k}\|_{L^2(\R^2)}
+
\|D g\|_{L^2(\R^2)}
\right).
\end{split}
\end{equation}
By combining~\eqref{eq:E_comp_prod} in \cref{res:E_comp} with the fact that $u_{s_k}=g$ on $A\setminus\Omega$, we find  that 
\begin{equation}
\label{strongout}
v_{s_k}\to g
\quad
\text{in}\
L^2(A\setminus\Omega)
\
\text{as}\
k\to\infty.
\end{equation}
Therefore, by exploiting~\eqref{strongout} in combination with~\eqref{eq:sup_F} and~\eqref{eq:E_comp_prod}, we get that
\begin{equation}
\label{eq:tozero_prod}
\lim_{k\to\infty}
\|D^{s_k} u_{s_k}\|_{L^2(\R^2)}\,\|v_{s_k}-g\|_{L^2(A\setminus\Omega)}
=0,
\end{equation}
which, together with~\eqref{eq:salire} and~\eqref{eq:mischia}, yields that $\|\mu\|_{\flute(A\setminus\Omega)}=0$ and thus $\supp \mu\subset\overline\Omega$.
Finally, we show that 
\begin{equation}
\label{eq:deg_bordo}
\mu(\overline\Omega)=\deg(g|_{\partial\Omega},\partial\Omega).
\end{equation}
To this end, let $\varphi\in\lip_c(A)$ be such that $\varphi=1$ on $\overline\Omega$.
According to \cref{def:j}, we can decompose
\begin{equation*}
j(v_{s_k})
=
j(v_{s_k},v_{s_k})
=
j(v_{s_k}-g,v_{s_k})
+
j(g,v_{s_k})
\end{equation*}
for $k\in\N$.
Since $\nabla\varphi=0$ on $\Omega$ by definition, from~\eqref{eq:tozero_prod} we deduce that
\begin{equation*}
\limsup_{k\to\infty}
\left|
\int_A 
j(v_{s_k}-g,v_{s_k})\times\nabla\varphi\di x
\right|
\le
\lip(\varphi)
\lim_{k\to\infty}
\|v_{s_k}-g\|_{L^2(A\setminus\Omega)}
\,
\|D v_{s_k}\|_{L^2(A)}
=0.
\end{equation*}
Thus, 
by combining this with~\eqref{eq:conv_int_jac}, \eqref{eq:cucu_jac} and~\eqref{eq:cucu_sym} in \cref{res:cucu}, and~\eqref{strongout}, we get
\begin{equation*}
\begin{split}
\pi\int_A\varphi\di\mu
&=
\lim_{k\to\infty}
\int_A\jac(v_{s_k})
\,
\varphi\di x
=
\lim_{k\to\infty}
\int_A 
j(v_{s_k})\times\nabla\varphi\di x
=
\lim_{k\to\infty}
\int_A 
j(g,v_{s_k})\times\nabla\varphi\di x
\\
&=
\lim_{k\to\infty}
\int_A 
j(v_{s_k},g)\times\nabla\varphi\di x
=
\int_A 
j(g)\times\nabla\varphi\di x
=
\pi
\int_A\jac(g)
\,
\varphi\di x.
\end{split}
\end{equation*} 
Since $\varphi=1$ on $\overline\Omega$, $\supp\mu\subset\overline\Omega$ and $|g|=1$ on $A\setminus\Omega$, by \cref{res:null_jac,res:deg} we hence get
\begin{equation*}
\mu(\overline\Omega)
=
\frac1{\pi}
\int_A\varphi\di\mu
=
\frac1{\pi}\int_A\jac(g)
\,
\varphi\di x
=
\frac1{\pi}\int_\Omega\jac(g)
\di x
=
\deg(g|_{\partial\Omega},\partial\Omega),
\end{equation*}
proving~\eqref{eq:deg_bordo} and thus yielding the conclusion.
\qed

\subsection{Proof of  claim \texorpdfstring{\ref{item:main_liminf}}{the liminf inequality} of \texorpdfstring{\cref{resi:main}}{the main result}}
\label{subsec:liminf}

Let $u_{s_k}\in H^{s_k}_g(\Omega;B_L)$, with $s_k\to1^-$ as $k\to\infty$, be such that $\jac^{s_k}(u_{s_k})\xflute{\overline\Omega}\pi\mu
$ as $s\to1^-$ for some $\mu\in\atoms(\overline\Omega)$.
Letting $v_{s_k}=I_{1-{s_k}}^Ru_{s_k}$ as in \cref{def:truncated_Riesz}, and repeating the argument of the proof of \cref{resi:main}\ref{item:main_comp}, we get 
\begin{equation*}
\mu(\overline\Omega)=\deg(\bar u|_{\partial\Omega},\partial\Omega)
\end{equation*}
and 
\begin{equation*}
\jac^{s_k}(v_{s_k})\xflute{A}\pi\nu
\quad\text{as}\
k\to\infty,
\end{equation*}
for some $\nu\in\atoms(A)$ such that $\supp \nu\subset\overline\Omega$ and $\nu|_{\overline\Omega}=\mu$.  
Therefore, owing to~\eqref{eq:E_comp_GL} in \cref{res:E_comp}, \eqref{eq:lim_gamma_R}, \eqref{eq:lim_Lambda} and \cref{res:GL}\ref{item:GL_liminf}, we can estimate
\begin{equation*}
\begin{split}
\liminf_{k\to\infty}
\GL_\lambda^{s_k}(u_{s_k};\Omega)
&\ge 
\liminf_{k\to\infty}
\frac{\eta}2 
\left(\frac{2\pi R^{1-{s_k}}}{(1-{s_k})\gamma_{s_k}}\right)^2
\,
\GL_{\sqrt{1-{s_k}},\Lambda_{{s_k},\lambda,\eta}}(v_{s_k};A)
\\
&=
\frac\eta2\,
\liminf_{k\to\infty}
\GL_{\sqrt{1-{s_k}},\Lambda_{{s_k},\eta,\lambda}}(v_{s_k};A)
\\
&\ge 
\frac\eta2\,
\liminf_{k\to\infty}
\GL_{\sqrt{1-{s_k}},\frac{\Lambda_{1,\lambda,\eta}}{2}}(v_{s_k};A)
\\
&\ge 
\eta\,\pi|\nu|(A)
=
\eta\,\pi|\mu|(\overline\Omega),
\end{split}
\end{equation*}
for any $\eta\in(0,1)$, yielding the conclusion.
\qed

\subsection{Truncated vortex}
\label{subsec:vortex}

The proof of point~\ref{item:main_limsup} of \cref{resi:main} requires some preliminaries.
We begin by introducing the following definition.

\begin{definition}[Cut-off function]
\label{def:cutoff}
Given $0\le r<R<\infty$, we define the \emph{cut-off function} $\eta_{r,R}\colon\R^2\to[0,1]$ by letting
\begin{equation}
\label{eq:cutoff}
\eta_{r,R}(t)
=
\begin{cases}
0
&
\text{if}\
|x|\in[0,r],\\[2ex]
\dfrac{|x|-r}{R-r}
&
\text{if}\
|x|\in[r,R],\\[2ex]
1
&
\text{if}\
|x|\ge R.
\end{cases}
\end{equation}
\end{definition}

In the following result, we collect some properties of the cut-off function $\eta_{r,R}$ given in \cref{def:cutoff} that will be useful below.

\begin{lemma}
\label{res:cutoff}
Given $0<r<R<\infty$, the function $\eta_{r,R}$ in \cref{def:cutoff} satisfies 
\begin{equation}
\label{eq:cutoff_lip}
|\eta_{r,R}|\le 1,
\quad 
\lip(\eta_{r,R})=\frac1{R-r},
\end{equation}
\begin{equation}
\label{eq:cutoff_ball}
\iint_{(B_R\times B_R)\cap\bigtriangleup_\tau}
\frac{|\eta_{r,R}(x)-\eta_{r,R}(y)|^2}{|x-y|^{2+2s}}\di x\di y
\le 
\frac{\pi^2 R^2}{(R-r)^2}
\,
\frac{\tau^{2-2s}}{1-s},
\end{equation}
\begin{equation}
\label{eq:cutoff_ring}
\iint_{(B_\rho\setminus B_r)\times(B_\rho\setminus B_r)\cap\bigtriangleup_\tau}
\frac{|\eta_{r,R}(x)-\eta_{r,R}(y)|^2}{|x-y|^{2+2s}}\di x\di y
\le
\frac{\pi^2((R+\tau)^2-r^2)}{(R-r)^2}
\,
\frac{\tau^{2-2s}}{1-s}
\end{equation}
and
\begin{equation}
\label{eq:cutoff_pozzo}
\int_{B_\varrho}\left(\eta_{r,R}^2-1\right)^2\di x
\le 
\pi R^2
\end{equation}
for $s\in(0,1)$, $\tau\in(0,r)$ and $\rho>R+\tau$.
\end{lemma}

In the proof of \cref{res:cutoff}, we will need the following elementary estimates.

\begin{lemma}
\label{res:tir}
Let $n\in\N$.
If $x,y\in\R^n\setminus\set*{0}$, then
\begin{equation}
\label{eq:tir}
\left|
\frac{x}{|x|}
-
\frac{y}{|y|}
\right|
\le 
\frac{2}{|x|}\,
|x-y|
\end{equation}
and
\begin{equation}
\label{eq:tirpro}
\left|
\frac{x}{|x|}
-
\frac{y}{|y|}
\right|^2
\le 
\frac{|x|}{|y|}
\left(
\frac{|y-x|^2}{|x|^2}
-
\frac{|x\cdot(y-x)|^2}{|x|^4}
+
\frac{|y-x|^3}{|x|^3}
\right).
\end{equation}
\end{lemma}

\begin{proof}
By the triangular inequality, we can estimate
\begin{equation*}
\left|
\frac{x}{|x|}
-
\frac{y}{|y|}
\right|
\le 
\left|
\frac{x}{|x|}
-
\frac{y}{|x|}
\right|
+
\left|
\frac{y}{|x|}
-
\frac{y}{|y|}
\right|
\le 
\frac{2}{|x|}
\,
|x-y|,
\end{equation*}
proving~\eqref{eq:tir}.
To prove~\eqref{eq:tirpro}, instead, we observe that
\begin{equation*}
\begin{split}
\left|
\frac{x}{|x|}
-
\frac{y}{|y|}
\right|^2
&=
2\left(1-\frac{x\cdot y}{|x|\,|y|}\right)
=
2\left(1-\frac{x\cdot (y-x)}{|x|\,|y|}
-
\frac{|x|}{|y|}
\right)
=
2\,\frac{|x|}{|y|}\left(\frac{|y|}{|x|}-\frac{x\cdot (y-x)}{|x|^2}
-
1
\right)
\\
&=
2\,\frac{|x|}{|y|}
\bigg(\sqrt{\left|\frac{x}{|x|}+\frac{y-x}{|x|}\right|^2}-1-\frac{x\cdot (y-x)}{|x|^2}\bigg)
\\
&=
2\,\frac{|x|}{|y|}
\bigg(\sqrt{1+2\,\frac{x\cdot(y-x)}{|x|^2}+\frac{|y-x|^2}{|x|^2}}-1-\frac{x\cdot (y-x)}{|x|^2}\bigg).
\end{split}
\end{equation*} 
Thanks to the elementary inequality
\begin{equation*}
\sqrt{1+t}\le 1+\frac12\,t-\frac18\,t^2
\quad
\text{for}\
t\ge0,
\end{equation*}
we thus get that
\begin{equation*}
\begin{split}
\left|
\frac{x}{|x|}
-
\frac{y}{|y|}
\right|^2
&\le
2\,\frac{|x|}{|y|}
\bigg(\frac12\,\frac{|y-x|^2}{|x|^2}
-\frac18
\left(
2\,\frac{x\cdot(y-x)}{|x|^2}+\frac{|y-x|^2}{|x|^2}
\right)^2
\bigg)
\\
&\le
2\,\frac{|x|}{|y|}
\left(\frac12\,\frac{|y-x|^2}{|x|^2}
-
\frac12\,
\frac{|x\cdot(y-x)|^2}{|x|^2}
-
\frac12\,
\frac{(x\cdot(y-x))\,|y-x|^2}{|x|^4}
\right)
\\
&\le 
\frac{|x|}{|y|}
\left(\frac{|y-x|^2}{|x|^2}
-
\frac{|x\cdot(y-x)|^2}{|x|^2}
+
\frac{|y-x|^3}{|x|^3}
\right),
\end{split}
\end{equation*}
yielding~\eqref{eq:tirpro} and concluding the proof.
\end{proof}

\begin{proof}[Proof of \cref{res:cutoff}]
From the definition in~\eqref{eq:cutoff}, we get~\eqref{eq:cutoff_lip} and thus
\begin{equation*}
\int_{B_\varrho}\left(\eta_{r,R}^2-1\right)^2\di x
=
\int_{B_R}\left(\eta_{r,R}^2-1\right)^2\di x
\le 
|B_R|,
\end{equation*}
giving~\eqref{eq:cutoff_pozzo}.
To prove~\eqref{eq:cutoff_ball}, we observe that
\begin{equation*}
\begin{split}
\int_{(B_R\times B_R)\cap\bigtriangleup_\tau}
&
\frac{|\eta_{r,R}(|x|)-\eta_{r,R}(|y|)|^2}{|x-y|^{2+2s}}\di x\di y
\le 
\int_{B_R}
\int_{B_\tau}
\frac{|\eta_{r,R}(|x|)-\eta_{r,R}(|x+h|)|^2}{|h|^{2+2s}}\di h\di x
\\
&
\le
\lip(\eta_{r,R})^2
\int_{B_R}
\int_{B_\tau}
\frac{|h|^2}{|h|^{2+2s}}\di h\di x
=
\frac{\pi^2 R^2}{(R-r)^2}
\,
\frac{\tau^{2-2s}}{1-s},
\end{split}
\end{equation*}
while, to prove~\eqref{eq:cutoff_ring}, we note that 
\begin{equation*}
x\in B_\varrho\setminus B_{R+\tau}
\implies
|x+h|>R
\quad
\text{for}\ h\in B_\tau,
\end{equation*}
and thus we can estimate
\begin{equation*}
\begin{split}
&\iint_{(B_\rho\setminus B_r)\times(B_\rho\setminus B_r)\cap\bigtriangleup_\tau}
\frac{|\eta_{r,R}(x)-\eta_{r,R}(y)|^2}{|x-y|^{2+2s}}\di x\di y
\\
&\le 
\int_{B_\tau}\frac{1}{|h|^{2+2s}}
\int_{B_\varrho\setminus B_r}
|\eta_{r,R}(x)-\eta_{r,R}(x+h)|^2
\di x\di h
\\
&=
\int_{B_\tau}
\frac{1}{|h|^{2+2s}}
\int_{B_{R+\tau}\setminus B_r}
|\eta_{r,R}(|x|)-\eta_{r,R}(|x+h|)|^2
\di x
\di h
\le 
\frac{\pi^2((R+\tau)^2-r^2)}{(R-r)^2}
\,
\frac{\tau^{2-2s}}{1-s},
\end{split}
\end{equation*}
and the proof is complete.
\end{proof}

We can introduce the notion of \emph{truncated vortex} by exploiting the cut-off function $\eta_{r,R}$ given in \cref{def:cutoff}. 

\begin{definition}[Truncated vortex]
\label{def:vortex}
Given $0\le r<R<\infty$ and $d\in\set*{-1,1}$, we let $\upsilon_{d,r,R}\colon\R^2\to\R^2$ be the \emph{truncated vortex} defined as
\begin{equation}
\label{eq:vortex}
\upsilon_{d,r,R}(x)
=
\eta_{r,R}(|x|)\,
\frac{x^d}{|x|},
\quad
\text{for}\
x\in\R^2,
\end{equation}
where $\eta_{r,R}$ is as in \cref{def:cutoff} and $x^d=\svec{x_1}{(-1)^dx_2}$ for all $x\in\R^2$.
\end{definition}

We now collect several properties of the truncated vortex introduced in \cref{def:vortex}.
On the one hand, we have the following result.

\begin{lemma}
\label{res:vortex1}
Given $0\le r<R<\infty$ and $d\in\set*{-1,1}$, the truncated vortex $\upsilon_{d,r,R}$ in \cref{def:vortex} satisfies
\begin{equation}
\label{eq:vortex_norm}
\upsilon_{d,r,R}
\in 
H^1_{\rm loc}(\R^2;\R^2)\cap L^\infty(\R^2;\R^2),
\end{equation}
with
\begin{equation}
\label{eq:vortex_grad_est}
|D\upsilon_{d,r,R}(x)|
\le 
\frac{C}{|x|}\,\eta_{r,R}(|x|),
\quad
x\in\R^2,
\end{equation}
where $C>0$ is a numerical constant.
In addition, it holds that
\begin{equation}
\label{eq:vortex_L2_dist}
\int_{\R^2}|\upsilon_{d,r,R}-\upsilon_{d,0,R}|^2\di x
\le4\pi R^2, 
\end{equation}
\begin{equation}
\label{eq:vortex_jac}
\jac(\upsilon_{d,0,R})
=
\pi d\,\frac{\chi_{B_R}}{|B_R|}
\end{equation} 
and 
\begin{equation}
\label{eq:vortex_H1}
\int_{B_\varrho}|D\upsilon_{d,0,R}|^2\di x
\le 
C\log\left(\frac{\varrho}{R}\right)
\end{equation}
for $\varrho>R$, where $C>0$ is a numerical constant.
\end{lemma}

\begin{proof}
From~\eqref{eq:vortex} in \cref{def:vortex}, we get that $\|\upsilon_{d,r,R}\|_{L^\infty}\le 1$ and a direct computation yields~\eqref{eq:vortex_grad_est}, from which we get~\eqref{eq:vortex_norm}.
The validity of~\eqref{eq:vortex_L2_dist} is a consequence of~\eqref{eq:cutoff} in \cref{def:cutoff}.
Finally, from the fact that 
\begin{equation*}
D\upsilon_{d,0,R}(x)
=
\frac{\chi_{B_R}}{R}\,J_d
+
\frac{\chi_{B_R^c}}{|x|}
\Big(J_d-\frac{x\otimes x^d}{|x|^2}\Big),
\quad
x\in\R^2\setminus\set*{0},
\end{equation*} 
where we have set
\begin{equation*}
J_d=\mat{1}{0}{0}{(-1)^d},
\quad
d\in\set*{1,-1},
\end{equation*}
we infer~\eqref{eq:vortex_jac} and~\eqref{eq:vortex_H1}, concluding the proof.
\end{proof}

On the other hand, similarly to \cref{res:cutoff}, we have the following result.

\begin{lemma}
\label{res:vortex2} Given $0<r<R<\infty$ and $d\in\set*{-1,1}$, the truncated vortex $\upsilon_{d,r,R}$ in \cref{def:vortex} satisfies
\begin{equation}
\label{eq:vortex_ball}
\iint_{(B_R\times B_R)\cap\bigtriangleup_\tau}
\frac{|\upsilon_{d,r,R}(x)-\upsilon_{d,r,R}(y)|^2}{|x-y|^{2+2s}}\di x\di y
\le
\Big(
\frac{8\pi^2(R^2-r^2)}{r^2}
+
\frac{2\pi^2 R^2}{(R-r)^2}
\Big)
\frac{\tau^{2-2s}}{1-s},
\end{equation}
\begin{equation}
\label{eq:vortex_ring}
\begin{split}
\iint_{((B_\varrho\setminus B_r)\times(B_\varrho\setminus B_r))\cap\bigtriangleup_\tau}
&
\frac{|\upsilon_{d,r,R}(x)-\upsilon_{d,r,R}(y)|^2}{|x-y|^{2+2s}}\di x\di y
\\
&\le
(1+\e)
\Big(
\frac{\pi^2\,r}{r-\tau}
\log\Big(\frac{\varrho}{r}\Big)
\frac{\tau^{2-2s}}{1-s}
+
\frac{4\pi^2(\varrho-r)}{\varrho\,(r-\tau)}
\,
\frac{\tau^{3-2s}}{3-2s}
\Big)
\\
&\quad+
\left(1+\frac1\e\right)
\frac{\pi^2((R+\tau)^2-r^2)}{(R-r)^2}
\,
\frac{\tau^{2-2s}}{1-s}
\end{split}
\end{equation}
and 
\begin{equation}
\label{eq:vortex_pozzo}
\int_{B_\varrho}\left(|\upsilon_{d,r,R}|^2-1\right)^2
\di x
\le 
\pi R^2
\end{equation}
for $s\in(0,1)$, $\tau\in(0,r)$, $\varrho>R+\tau$ and $\e>0$.
\end{lemma}

\begin{proof}
We have that 
\begin{equation}
\label{alice}
|\upsilon_{d,r,R}(x)-\upsilon_{d,r,R}(y)|
\le 
\eta_{r,R}(|x|)
\left|
\frac{x}{|x|}
-
\frac{y}{|y|}
\right|
+
|\eta_{r,R}(|x|)-\eta_{r,R}(|y|)|.
\end{equation}
Thanks to~\eqref{eq:tir} in \cref{res:tir}, we can estimate
\begin{equation}
\label{panco}
\begin{split}
\iint_{(B_R\times B_R)\cap\bigtriangleup_\tau}
&
\eta_{r,R}^2(|x|)
\,\frac{
\left|
\frac{x}{|x|}
-
\frac{y}{|y|}
\right|^2}{|x-y|^{2+2s}}\di x\di y
\le 
\int_{((B_R\setminus B_r)\times B_R)\cap\bigtriangleup_\tau}
\frac{4}{|x|^2}
\,
\frac{|x-y|^2}{|x-y|^{2+2s}}\di x\di y
\\
&
\le
\frac{4}{r^2}
\int_{B_R\setminus B_r}
\int_{B_\tau}
\frac{|h|^2}{|h|^{2+2s}}
\di h
\di x 
=
\frac{4\pi^2(R^2-r^2)}{r^2}
\,
\frac{\tau^{2-2s}}{1-s}.
\end{split}
\end{equation}
Therefore, by combining~\eqref{alice} and~\eqref{panco} with~\eqref{eq:cutoff_ball} in \cref{res:cutoff}, we get that
\begin{equation*}
\begin{split}
&\iint_{(B_R\times B_R)\cap\bigtriangleup_\tau}
\frac{|\upsilon_{d,r,R}(x)-\upsilon_{d,r,R}(y)|^2}{|x-y|^{2+2s}}\di x\di y
%\\
%&\le
%2
%\int_{(B_R\times B_R)\cap\bigtriangleup_\tau}
%\eta_{r,R}^2(|x|)
%\,\frac{
%\left|
%\frac{x}{|x|}
%-
%\frac{y}{|y|}
%\right|^2}{|x-y|^{2+2s}}\di x\di y
%+
%2 
%\int_{(B_R\times B_R)\cap\bigtriangleup_\tau}
%\frac{|\eta_{r,R}(|x|)-\eta_{r,R}(|y|)|^2}{|x-y|^{2+2s}}\di x\di y
%\\
%&
\le
\left(
\frac{8\pi^2(R^2-r^2)}{r^2}
+
\frac{2\pi^2 R^2}{(R-r)^2}
\right)
\frac{\tau^{2-2s}}{1-s},
\end{split}
\end{equation*}
proving~\eqref{eq:vortex_ball}.
We now deal with~\eqref{eq:vortex_ring}.
In view of the elementary inequality
\begin{equation*}
(a+b)^2\le(1+\e)\,a^2+\left(1+\frac1{\e}\right)\,b^2
\quad
\text{for}\
a,b\ge0,\ \e>0,
\end{equation*}
we can revisit~\eqref{alice} as
\begin{equation}
\label{rici}
\begin{split}
|\upsilon_{d,r,R}(x)-\upsilon_{d,r,R}(y)|^2
&\le 
(1+\e)
\,
\eta_{r,R}(|x|)^2
\left|
\frac{x}{|x|}
-
\frac{y}{|y|}
\right|^2
+
\left(1+\frac1\e\right)
|\eta_{r,R}(|x|)-\eta_{r,R}(|y|)|^2
\end{split}
\end{equation}
whenever $\e>0$.
We now observe that
\begin{equation}
\label{clag}
\begin{split}
\iint_{((B_\varrho\setminus B_r)\times(B_\varrho\setminus B_r))\cap\bigtriangleup_\tau}
&
\eta_{r,R}^2(|x|)
\,\frac{
\left|
\frac{x}{|x|}
-
\frac{y}{|y|}
\right|^2}{|x-y|^{2+2s}}\di x\di y
\\
&\le
\int_{B_\tau}
\frac{1}{|h|^{2+2s}} 
\int_{B_\varrho\setminus B_r}
\left|\frac{x}{|x|}-\frac{x+h}{|x+h|}\right|^2
\di x\di h.
\end{split}
\end{equation}
Since $\tau<r$, we have that 
\begin{equation*}
|x+h|
\ge
|x|-|h|
\ge 
r-\tau>0
\quad
\text{for}\
x\in B_\varrho\setminus B_r,\ h\in B_\tau 
\end{equation*}
and therefore we can estimate
\begin{equation*}
\frac{|x|}{|x+h|}
=
\frac{1}{\left|\frac{x}{|x|}+\frac{h}{|x|}\right|}
\le 
\frac{1}{1-\frac{|h|}{|x|}}
\le 
\frac{1}{1-\frac{\tau}{r}}
=
\frac{r}{r-\tau}
\quad
\text{for}\
x\in B_\varrho\setminus B_r,\ h\in B_\tau.
\end{equation*}
Thus, thanks to~\eqref{eq:tirpro} in \cref{res:tir}, we get that
\begin{equation*}
\begin{split}
\int_{B_\tau}
&
\frac{1}{|h|^{2+2s}}
\int_{B_\varrho\setminus B_r}
\left|
\frac{x}{|x|}
-
\frac{x+h}{|x+h|}
\right|^2
\di x\di h
\\
&\le
\int_{B_\tau}
\frac{1}{|h|^{2+2s}}
\int_{B_\varrho\setminus B_r}
\frac{|x|}{|x+h|}
\left(
\frac{|h|^2}{|x|^2}
-
\frac{|x\cdot h|^2}{|x|^4}+\frac{|h|^3}{|x|^3}
\right)
\di x\di h
\\
&\le
\frac{r}{r-\tau}
\int_{B_\tau}
\frac{1}{|h|^{2+2s}}
\int_{B_\varrho\setminus B_r}
\frac{|h|^2}{|x|^2}
-
\frac{|x\cdot h|^2}{|x|^4}+\frac{|h|^3}{|x|^3}
\di x\di h.
\end{split}
\end{equation*}
At this point, we make the following observation. 
Given $h\in B_\tau$, we can find a rotation matrix $\mathcal R\in\mathrm{SO}(2)$ such that $h=\mathcal R\mathrm e_1$. 
Therefore, we can compute
\begin{equation*}
\begin{split}
\int_{B_\varrho\setminus B_r}
\frac{|x\cdot h|^2}{|x|^4}\di x
=
\int_{B_\varrho\setminus B_r}
\frac{|x\cdot \mathrm e_1|^2}{|x|^4}\di x
=
\int_{B_\varrho\setminus B_r}
\frac{x_1^2}{|x|^4}\di x
\end{split}
\end{equation*}
and thus, by replacing $\mathrm e_1$ with $\mathrm e_2$, we get that
\begin{equation*}
\int_{B_\varrho\setminus B_r}
\frac{|x\cdot h|^2}{|x|^4}\di x
=
\frac12
\int_{B_\varrho\setminus B_r}
\frac{x_1^2}{|x|^4}\di x
+
\frac12
\int_{B_\varrho\setminus B_r}
\frac{x_2^2}{|x|^4}\di x
=
\frac12
\int_{B_\varrho\setminus B_r}
\frac{|x|^2}{|x|^4}\di x
=
\pi\,\log\left(\frac{\varrho}{r}\right).
\end{equation*}
As a consequence, we have that
\begin{equation*}
\int_{B_\tau}
\frac{1}{|h|^{2+2s}}
\int_{B_\varrho\setminus B_r}
\frac{|h|^2}{|x|^2}
-
\frac{|x\cdot h|^2}{|x|^4}+\frac{|h|^3}{|x|^3}
\di x\di h
=
\pi^2
\log\left(\frac{\varrho}{r}\right)
\,
\frac{\tau^{2-2s}}{1-s}
+
4\pi^2
\left(\frac1r-\frac1\varrho\right)
\frac{\tau^{3-2s}}{3-2s},
\end{equation*}
from which we deduce that 
\begin{equation}
\label{gio}
\int_{B_\tau}
\frac{1}{|h|^{2+2s}}
\int_{B_\varrho\setminus B_r}
\left|
\frac{x}{|x|}
-
\frac{x+h}{|x+h|}
\right|^2
\di x\di h
\le 
\frac{\pi^2\,r}{r-\tau}
\log\left(\frac{\varrho}{r}\right)
\frac{\tau^{2-2s}}{1-s}
+
\frac{4\pi^2(\varrho-r)}{\varrho\,(r-\tau)}
\,
\frac{\tau^{3-2s}}{3-2s}.
\end{equation}
Thus, by combining~\eqref{rici}, \eqref{clag} and~\eqref{gio} with~\eqref{eq:cutoff_ring} in \cref{res:vortex2}, we get that
\begin{equation*}
\begin{split}
\iint_{((B_\varrho\setminus B_r)\times(B_\varrho\setminus B_r))\cap\bigtriangleup_\tau}
&
\frac{|\upsilon_{d,r,R}(x)-\upsilon_{d,r,R}(y)|^2}{|x-y|^{2+2s}}\di x\di y
\\
&\le
(1+\e)
\Big(
\frac{\pi^2\,r}{r-\tau}
\log\left(\frac{\varrho}{r}\right)
\frac{\tau^{2-2s}}{1-s}
+
\frac{4\pi^2(\varrho-r)}{\varrho\,(r-\tau)}
\,
\frac{\tau^{3-2s}}{3-2s}
\Big)
\\
&\quad+
\left(1+\frac1\e\right)
\frac{\pi^2((R+\tau)^2-r^2)}{(R-r)^2}
\,
\frac{\tau^{2-2s}}{1-s}
\end{split}
\end{equation*}
whenever $\e>0$, yielding~\eqref{eq:vortex_ring}.
The validity of~\eqref{eq:vortex_pozzo} follows by combining the definition in~\eqref{eq:vortex} with~\eqref{eq:cutoff_pozzo}.
The proof is complete.
\end{proof}

\subsection{Proof of claim \texorpdfstring{\ref{item:main_limsup}}{the limsup} of \texorpdfstring{\cref{resi:main}}{the main result}}
\label{subsec:limsup}

Below, in order to avoid heavy notation, we will frequently adopt the following shorthand.
Given  $x_1,\dots,x_N\in\R^2$ for some $N\in\N$, any $\ell>0$ and any non-empty open set $V\subset\R^2$, we let 
\begin{equation}
\label{eq:bucato}
\widehat V_\ell
=
V\setminus\Big(\bigcup_{i=1}^N \overline{B_\ell(x_i)}\Big).
\end{equation}

\begin{proof}[Proof of claim {\rm\ref{item:main_limsup}} of \cref{resi:main}]
We divide the proof in three steps.

\vspace{1ex}

\textit{Step 1}.
Let $\mu\in\atoms(\overline\Omega)$ be given by
\begin{equation}
\label{eq:special}
\mu=\sum_{i=1}^N d_i\,\delta_{x_i}
\end{equation}
for some $N\in\N$, $x_i\in\Omega$ and $d_i\in\set*{-1,1}$, for each $i\in\set*{1,\dots,N}$, such that 
\begin{equation}
\label{eq:special_deg}
\mu(\overline\Omega)
=
\mu(\Omega)
=
\sum_{i=1}^Nd_i=\deg(g|_{\partial\Omega},\partial\Omega).
\end{equation}
In this step, we construct $u_s\in H^s_g(\Omega;B_L)$ for $s\in(0,1)$ such that 
\begin{equation}
\label{eq:special_limsup}
\limsup_{s\to1^-}
\GL^s_\lambda(u_s)
\le 
\pi|\mu|(\overline\Omega)
=
\pi N
\end{equation}
whenever $\lambda\in(0,\infty)$.

To define $u_s$, we need to introduce some parameters.
We let
\begin{equation}
\label{eq:tau_M}
\tau_s=\sqrt{1-s}
\quad
\text{and}
\quad
M_s=\sqrt[4]{|\log(1-s)|}
\quad
\text{for}\
s\in(0,1).
\end{equation}
In view of the assumption made in~\eqref{eq:special}, we define 
\begin{equation}
\label{eq:spazio}
\overline r
=
\min\set*{\operatorname{dist}(x_i,\R^2\setminus\Omega),\ \frac12\,|x_i-x_j| : i,j\in\set*{1,\dots,N},\ i\ne j}\in(0,\infty)
\end{equation}
and we fix
\begin{equation}
\label{eq:raggio}
r<\frac{\overline r}2.
\end{equation}
Without loss of generality, we tacitly work with $s\in(0,1)$ sufficiently close to~$1$ so that
\begin{equation}
\label{eq:piccolo}
(M_s+2)\,\tau_s<\frac r2,
\end{equation}
which is always possibile in virtue of the definitions in~\eqref{eq:tau_M}.

We can now define $u_s$.
With the notation of \cref{def:vortex}, we define 
\begin{equation}
\label{eq:recovery_i}
u_{s,i}
=
\upsilon_{d_i,M_s\tau_s,(M_s+1)\tau_s}(\,\cdot\,-x_i)
\quad
\text{for each}\ i\in\set*{1,\dots,N}.
\end{equation}
By well-known results (e.g., see~\cite{BBH94}*{Th.~I.4}), due to~\eqref{eq:special_deg},
we can find $\widehat u\in H^1(\widehat\Omega_r;\sphere)$ (where $\widehat\Omega_r$ is defined using the shorthand~\eqref{eq:bucato} with points $x_1,\dots,x_N$ given by~\eqref{eq:special} and $\ell=r$ as fixed in~\eqref{eq:raggio}) such that
\begin{equation}
\label{eq:recovery_hat}
\widehat u=
\begin{cases}
u_{s,i}
&
\text{on}\
\partial B_r(x_i),
\
\text{for each}\
i\in\set*{1,\dots,N},
\\[1ex]
g
& 
\text{on}\ \partial\Omega.
\end{cases}
\end{equation}
Since, by~\eqref{eq:recovery_i}  and by \cref{def:truncated_Riesz}, for each $i\in\set*{1,\dots,N}$, we have 
\begin{equation*}
u_{s,i}=\frac{(\,\cdot\,-x_i)^{d_i}}{|\,\cdot\,-x_i|}
\quad
\text{in an open neighborhood of}\ 
\partial B_r(x_i),
\end{equation*}
the function~$\widehat u$ does not depend on $s$. 
We thus define $u_s\colon\R^2\to\R^2$ by letting 
\begin{equation}
\label{eq:recovery}
u_s
=
\begin{cases}
u_{s,i} 
&
\text{in}\
\overline{B_r(x_i)},\ \text{for each}\ i\in\set*{1,\dots,N},
\\[1ex]
\widehat u
&
\text{in}\
\widehat\Omega_r,
\\[1ex]
g
&
\text{in}\ \R^2\setminus\Omega.
\end{cases}
\end{equation}

We observe that $u_s\in H^s_g(\Omega;B_L)$, since, by~\eqref{eq:recovery_i} and by~\eqref{eq:vortex_norm} in \cref{res:vortex1}, and by~\eqref{eq:recovery_hat}, we have that $u_s\in H^1(\R^2;\R^2)$ with $\|u_s\|_{L^\infty}\le\|g\|_{L^\infty}$.

We now detail the proof of~\eqref{eq:special_limsup}.
We begin by observing that, owing to~\eqref{eq:recovery}, the fact that $|u_s|=|\widehat u|=1$ on $\widehat\Omega_r$ and $|u_s|=|u_{s,i}|=1$ on $B_r(x_i)\setminus B_{(M_s+1)\tau_s}(x_i)$ for each $i\in\set*{1,\dots,N}$,  and to~\eqref{eq:vortex_pozzo} in \cref{res:vortex2}, we can estimate   
\begin{equation*}
\begin{split}
\int_\Omega\left(|u_s|^2-1\right)^2\di x
&=
\int_{\widehat\Omega_r}\left(|\widehat u|^2-1\right)^2\di x
+
\sum_{i=1}^N
\int_{B_r(x_i)}\left(|u_{s,i}|^2-1\right)^2\di x
\\
&=
N
\int_{B_{(M_s+1)\tau_s}}\left(|\upsilon_{d_i,M_s\tau_s,(M_s+1)\tau_s}|^2-1\right)^2\di x
\le
N\pi\,(M_s+1)^2\,\tau^2_s,
\end{split}
\end{equation*}
from which, in virtue of the definitions in~\eqref{eq:tau_M}, we get that
\begin{equation*}
\limsup_{s\to1^-}
\int_\Omega (|u_s|^2-1)^2\di x
=0.
\end{equation*}
Hence, to show~\eqref{eq:special_limsup}, by~\eqref{eq:relation_s-grad_seminorm} in \cref{res:equiv-frac}, we just need to prove that 
\begin{equation}
\label{eq:recovery_energy}
\limsup_{s\to1^-}
\frac{(1-s)}{|\log(1-s)|}\,[u_s]^2_s
\le 
\frac{\pi^2N}{2}.
\end{equation}
Recalling the notation introduced in~\eqref{eq:bermuda}, by \cref{res:fartau} and~\eqref{eq:recovery}, we can estimate
\begin{equation*}
\begin{split}
\iint_{(\R^2\times\R^2)\setminus \bigtriangleup_{\tau_s}}\frac{|u_s(x)-u_s(y)|^2}{|x-y|^{2+2s}}\di x \di y
&\le 
\frac{4\pi\|u_s\|_{L^2}^2}{s\,(1-s)^s}
\\
&\le 
\frac{4\pi}{s\,(1-s)^s}
\,
\left(\|g\|_{L^2}^2+\|\widehat u\|_{L^2(\widehat\Omega_r)}^2+N\pi r^2\right),
\end{split}
\end{equation*}
from which we deduce that 
\begin{equation}
\label{eq:recovery_fartau}
\lim_{s\to1^-}
\frac{(1-s)}{|\log(1-s)|}
\,
\iint_{(\R^2\times\R^2)\setminus \bigtriangleup_{\tau_s}}\frac{|u_s(x)-u_s(y)|^2}{|x-y|^{2+2s}}\di x \di y
=
0.
\end{equation}
Hence, by combining~\eqref{eq:recovery_energy} with~\eqref{eq:recovery_fartau}, the validity of~\eqref{eq:special_limsup} reduces to
\begin{equation}
\label{eq:recovery_tau}
\limsup_{s\to1^-}
\frac{(1-s)}{|\log(1-s)|}
\iint_{(\R^2\times\R^2)\cap \bigtriangleup_{\tau_s}}\frac{|u_s(x)-u_s(y)|^2}{|x-y|^{2+2s}}\di x\di y
\le 
\frac{\pi^2N}{2}.
\end{equation}

Our aim is now to estimate the integral 
\begin{equation*}
\iint_{\bigtriangleup_{\tau_s}}\frac{|u_s(x)-u_s(y)|^2}{|x-y|^{2+2s}}\di x\di y.
\end{equation*}
To do so, we let
\begin{equation*}
Q_s(V,W)
=
\iint_{(V\times W)\cap\bigtriangleup_{\tau_s}}\frac{|u_s(x)-u_s(y)|^2}{|x-y|^{2+2s}}\,\di x\di y
\end{equation*}
for any two measurable sets $V,W\subset\R^2$. 
Since we can write $\R^2
=
V_{s,1}\cup V_{s,2}\cup V_{s,3}$, with pairwise disjoint union, where
\begin{equation*}
V_{s,1}
=
\bigcup_{i=1}^N B_{(M_s+1)\tau_s}(x_i)
,
\quad
V_{s,2}
=
\bigcup_{i=1}^N B_r(x_i)\setminus B_{(M_s+1)\tau_s}(x_i)
\quad\text{and}\quad
V_{s,3}=\widehat{(\R^2)}_r,
\end{equation*}
and $\widehat{(\R^2)}_r$ is as in~\eqref{eq:bucato} (with $V=\R^2$ and $\ell=r$), we can decompose
\begin{equation*}
\iint_{\bigtriangleup_{\tau_s}}\frac{|u_s(x)-u_s(y)|^2}{|x-y|^{2+2s}}\di x\di y
=
\sum_{j,k=1}^3
Q_s(V_{s,j},V_{s,k}).
\end{equation*}
We now deal with each possible pair.  
We begin by observing that, if $x\in B_{(M_s+1)\tau_s}(x_i)$ and $y\in B_{(M_s+1)\tau_s}(x_h)$ for some $i,h\in\set*{1,\dots,N}$ with $i\ne h$, then
\begin{equation*}
|x-y|
> 
|x_i-x_h|-2(M_s+1)\tau_s
\ge 
2\bar r-
2(M_s+1)\tau_s
>
2r-
2(M_s+1)\tau_s
>\tau_s,
\end{equation*}
because of~\eqref{eq:spazio}, \eqref{eq:raggio} and~\eqref{eq:piccolo}, so that 
\begin{equation*}
Q_s(V_{s,1},V_{s,1})
=
\sum_{i=1}^N
Q_s
\left(
B_{(M_s+1)\tau_s}(x_i)
, B_{(M_s+1)\tau_s}(x_i)
\right).
\end{equation*}
Next, we notice that
\begin{equation*}
Q_s(V_{s,1},V_{s,2})
=
Q_s(\widetilde V_{s,1},V_{s,2})
\quad\text{and}\quad
Q_s(V_{s,2},V_{s,1})
=
Q_s(V_{s,2},\widetilde V_{s,1}),
\end{equation*}
where we set
\begin{equation*}
\widetilde V_{s,1}
=
\bigcup_{i=1}^N B_{(M_s+1)\tau_s}(x_i)
\setminus
B_{M_s\tau_s}(x_i).
\end{equation*}
Moreover, if $x\in B_r(x_i)$ and $y\in B_{(M_s+1)\tau_s}(x_h)$ for some $i,h\in\set*{1,\dots,N}$ with $i\ne h$, then 
\begin{equation*}
|x-y|
>
|x_i-x_h|-r-(M_s+1)\tau_s
\ge 
2\bar r-r-(M_s+1)\tau_s
>
r-(M_s+1)\tau_s>\tau_s,
\end{equation*}
again because of~\eqref{eq:spazio}, \eqref{eq:raggio} and~\eqref{eq:piccolo}, so that 
\begin{equation*}
Q_s(V_{s,1},V_{s,2})
=
\sum_{i=1}^N
Q_s
\left(
B_{(M_s+1)\tau_s}(x_i)
\setminus
B_{M_s\tau_s}(x_i)
,
B_r(x_i)
\setminus
B_{(M_s+1)\tau_s}(x_i)
\right)
\end{equation*}
and, similarly,
\begin{equation*}
Q_s(V_{s,2},V_{s,1})
=
\sum_{i=1}^N
Q_s
\left(
B_r(x_i)
\setminus
B_{(M_s+1)\tau_s}(x_i)
,
B_{(M_s+1)\tau_s}(x_i)
\setminus
B_{M_s\tau_s}(x_i)
\right).
\end{equation*} 
As above, if $x\in B_r(x_i)$ and $y\in B_r(x_h)$ for some $i,h\in\set*{1,\dots,N}$ with $i\ne h$, then
\begin{equation*}
|x-y|
>
|x_i-x_h|-2r
\ge 
2\bar r-2r
>
4r-2r
>
\tau_s,
\end{equation*}
again because of~\eqref{eq:spazio}, \eqref{eq:raggio} and~\eqref{eq:piccolo}, so that 
\begin{equation*}
Q_s(V_{s,2},V_{s,2})
=
\sum_{i=1}^N
Q_s
\left(
B_r(x_i)
\setminus
B_{(M_s+1)\tau_s}(x_i)
,
B_r(x_i)
\setminus
B_{(M_s+1)\tau_s}(x_i)
\right).
\end{equation*} 
We hence infer that 
\begin{equation*}
\begin{split}
Q_s&(V_{s,1},V_{s,2})
+
Q_s(V_{s,2},V_{s,1})
+
Q_s(V_{s,2},V_{s,2})
\\
&=
\sum_{i=1}^N
Q_s
\left(
B_{(M_s+1)\tau_s}(x_i)
\setminus
B_{M_s\tau_s}(x_i)
,
B_r(x_i)
\setminus
B_{(M_s+1)\tau_s}(x_i)
\right)
\\
&\quad+
\sum_{i=1}^N
Q_s
\left(
B_r(x_i)
\setminus
B_{(M_s+1)\tau_s}(x_i)
,
B_{(M_s+1)\tau_s}(x_i)
\setminus
B_{M_s\tau_s}(x_i)
\right)
\\
&\quad+
\sum_{i=1}^N
Q_s
\left(
B_r(x_i)
\setminus
B_{(M_s+1)\tau_s}(x_i)
,
B_r(x_i)
\setminus
B_{(M_s+1)\tau_s}(x_i)
\right)
\\
&=
\sum_{i=1}^N
Q_s
\left(
B_{(M_s+1)\tau_s}(x_i)
\setminus
B_{M_s\tau_s}(x_i)
,
B_r(x_i)
\setminus
B_{(M_s+1)\tau_s}(x_i)
\right)
\\
&\quad+
\sum_{i=1}^N
Q_s
\left(
B_r(x_i)
\setminus
B_{(M_s+1)\tau_s}(x_i)
,
B_r(x_i)
\setminus
B_{M_s\tau_s}(x_i)
\right)
\\
&\le 
\sum_{i=1}^N
Q_s
\left(
B_r(x_i)
\setminus
B_{M_s\tau_s}(x_i)
,
B_r(x_i)
\setminus
B_{M_s\tau_s}(x_i)
\right).
\end{split}
\end{equation*}
Finally, if $x\in B_{(M_s+1)\tau_s}(x_i)$ and $y\in\R^2\setminus B_r(x_h)$ for some $i,h\in\set*{1,\dots,N}$, then
\begin{equation*}
|x-y|\ge r-(M_s+1)\tau_s>\tau_s
\end{equation*}  
by~\eqref{eq:piccolo}, so that
\begin{equation*}
Q_s(V_{s,1},V_{s,3})
=
Q_s(V_{s,3},V_{s,1})
=
0.
\end{equation*}
Moreover, if $x\in V_{2,s}$ and $y\in V_{3,s}$, then $(x,y)\in\bigtriangleup_{\tau_s}$ only in the case in which
\begin{equation*}
x\in B_r(x_i)\setminus B_{r-\tau_s}(x_i)
\quad
\text{and}
\quad
y\in\R^2\setminus B_r(x_i)
\end{equation*}
for some $i\in\set*{1,\dots,N}$.
We therefore get that 
\begin{equation*}
\begin{split}
Q_s
&
(V_{2,s},V_{3,s})
+
Q_s(V_{3,s},V_{2,s})
+
Q_s(V_{3,s},V_{3,s})
=
Q_s
\left(
\bigcup_{i=1}^N
B_r(x_i)\setminus B_{r-\tau_s}(x_i)
,
\widehat{(\R^2)}_r
\right)
\\
&\quad+
Q_s
\left(
\widehat{(\R^2)}_r
,
\bigcup_{i=1}^N
B_r(x_i)\setminus B_{r-\tau_s}(x_i)
\right)
+
Q_s
\left(
\widehat{(\R^2)}_r
,
\widehat{(\R^2)}_r
\right)
\\
&=
Q_s
\left(
\bigcup_{i=1}^N
B_r(x_i)\setminus B_{r-\tau_s}(x_i)
,
\widehat{(\R^2)}_r
\right)
+
Q_s
\left(
\widehat{(\R^2)}_r
,
\widehat{(\R^2)}_{r-\tau_s}
\right)
\\
&\le 
Q_s
\left(
\widehat{(\R^2)}_{r-\tau_s}
,
\widehat{(\R^2)}_{r-\tau_s}
\right).
\end{split}
\end{equation*}
By combining all the above estimates, we thus conclude that 
\begin{equation}
\label{eq:pezzi}
\begin{split}
\iint_{\bigtriangleup_{\tau_s}}\frac{|u_s(x)-u_s(y)|^2}{|x-y|^{2+2s}}\di x\di y
&\le
\sum_{i=1}^N
Q_s
\left(
B_{(M_s+1)\tau_s}(x_i)
, B_{(M_s+1)\tau_s}(x_i)
\right)
\\
&\quad+
\sum_{i=1}^N
Q_s
\left(
B_r(x_i)
\setminus
B_{M_s\tau_s}(x_i)
,
B_r(x_i)
\setminus
B_{M_s\tau_s}(x_i)
\right)
\\
&\quad
+
Q_s
\big(
\widehat{(\R^2)}_{r-\tau_s}
,
\widehat{(\R^2)}_{r-\tau_s}
\big).
\end{split}
\end{equation} 

We estimate each term in the right-hand side of~\eqref{eq:pezzi} separately.
Recalling~\eqref{eq:recovery_i}, by~\eqref{eq:vortex_ball} in \cref{res:vortex2} we have that
\begin{equation}
\label{eq:pezzo1}
Q_s
\left(
B_{(M_s+1)\tau_s(x_i)}
, B_{(M_s+1)\tau_s(x_i)}
\right)
\le 
\frac{CM_s^2}{1-s}
\end{equation} 
for each $i\in\set*{1,\dots,N}$, 
where $C>0$ is a numerical constant.
Moreover, by~\eqref{eq:vortex_ring} in \cref{res:vortex2}, we similarly get that
\begin{equation}
\label{eq:pezzo2}
\begin{split}
Q_s
&
\left(
B_r(x_i)
\setminus
B_{M_s\tau_s}(x_i)
,
B_r(x_i)
\setminus
B_{M_s\tau_s}(x_i)
\right)
\\
&\le
(1+\e)
\,
\pi^2
\,
\frac{M_s\,
(1-s)^{1-s}}{M_s-1}
\,
\frac{\log r-\log M_s+\frac12|\log(1-s)|}{1-s}
+
\frac{CM_s}{\e\,(1-s)}
\end{split}
\end{equation}
for each $i\in\set*{1,\dots,N}$ and any $\e>0$,
where $C>0$ is a numerical constant. 
Finally, by the Fundamental Theorem of Calculus, Jensen's inequality and~\eqref{eq:piccolo}, we have that 
\begin{equation*}
\begin{split}
Q_s
\big(
\widehat{(\R^2)}_{r-\tau_s}
,
\widehat{(\R^2)}_{r-\tau_s}
\big)
&\le 
\int_{B_{\tau_s}}
\frac{1}{|h|^{2+2s}}
\int_{\widehat{(\R^2)}_{r-\tau_s}}
|u_s(x+h)-u_s(x)|^2
\di x
\di h
\\
&\le 
\int_{B_{\tau_s}}
\frac{|h|^2}{|h|^{2+2s}}
\int_{\widehat{(\R^2)}_{r-\tau_s}}
\left(\int_0^1|D u_s(x+th)|\di t
\right)^2
\di x
\di h
\\
&\le
\int_{B_{\tau_s}}
\frac{1}{|h|^{2s}}
\int_{\widehat{(\R^2)}_{r-\tau_s}}
\int_0^1|D u_s(x+th)|^2\di t
\di x
\di h
\\
&\le
\int_{B_{\tau_s}}
\frac{1}{|h|^{2s}}
\int_{\widehat{(\R^2)}_{r-2\tau_s}}
|D u_s(z)|^2
\di z
\di h
\\
&\le
\frac{C}{1-s}
\int_{\widehat{(\R^2)}_{r/2}}
|D u_s(z)|^2
\di z,
\end{split}
\end{equation*}
where $C>0$ is a numerical constant.
We also observe that, by~\eqref{eq:recovery}, we can write
\begin{equation*}
\begin{split}
\int_{\widehat{(\R^2)}_{r/2}}
|Du_s(z)|^2
\di z
&=
\int_{\R^2\setminus\Omega}
|D\bar u(z)|^2
\di z
+
\int_{\widehat\Omega_r}
|D\widehat u(z)|^2
\di z
+\sum_{i=1}^N
\int_{B_r(x_i)\setminus B_{r/2}(x_i)}|Du_{s,i}|^2\di x,
\end{split}
\end{equation*}
with, thanks to~\eqref{eq:vortex_grad_est} and again~\eqref{eq:piccolo},
\begin{equation*}
\int_{B_r(x_i)\setminus B_{r/2}(x_i)}
|Du_{s,i}|^2\di x
\le
C
\int_{B_r\setminus B_{r/2}}\frac{1}{|x|^2}\di x
\le 
4C\pi,
\end{equation*}
where $C>0$ is a numerical constant.
Therefore, we have that 
\begin{equation}
\label{eq:pezzo3}
Q_s
\big(
\widehat{(\R^2)}_{r-\tau_s}
,
\widehat{(\R^2)}_{r-\tau_s}
\big)
\le 
\frac{C}{1-s},
\end{equation}
where $C>0$ does not depend on~$s$.
Hence, by combining~\eqref{eq:pezzi}, \eqref{eq:pezzo1}, \eqref{eq:pezzo2} and~\eqref{eq:pezzo3}, and by recalling~\eqref{eq:tau_M}, we infer that
\begin{equation*}
\limsup_{s\to1^-}
\frac{(1-s)}{|\log(1-s)|}
\iint_{\bigtriangleup_{\tau_s}}\frac{|u_s(x)-u_s(y)|^2}{|x-y|^{2+2s}}\di x\di y
\le 
(1+\e)\,\frac{\pi^2 N}{2}
\end{equation*}
whenever $\e>0$, proving~\eqref{eq:recovery_tau} and thus~\eqref{eq:special_limsup}.

\vspace{1ex}

\textit{Step~2}.
Let $\mu$ be as in~\eqref{eq:special} and let $u_s$ be as in~\eqref{eq:recovery}. 
In this step, we prove that  
\begin{equation}
\label{eq:recovery_jac}
\jac^s(u_s)\xflute{\overline\Omega}\pi\mu
\quad
\text{as}\
s\to1^-.
\end{equation}
Thanks to~\eqref{eq:lim_gamma_R} and~\eqref{eq:jacs}, to show~\eqref{eq:recovery_jac} we just need to prove that 
\begin{equation}
\label{eq:recovery_I_jac}
\jac(v_s)\xflute{\overline\Omega}\pi\mu
\quad
\text{as}\
s\to1^-,
\end{equation}
where, following \cref{def:truncated_Riesz}, $v_s=I_{1-s}^R u_s$ for $s\in(0,1)$. 
To prove~\eqref{eq:recovery_I_jac}, in turn, we need to introduce the auxiliary function $\widetilde u_s\colon\R^2\to\R^2$ defined as
\begin{equation}
\label{eq:aux}
\widetilde u_s
=
\begin{cases}
\upsilon_{d_i,0,(M_s+1)\tau_s}(\,\cdot\,-x_i) 
&
\text{in}\
\overline{B_r(x_i)},\ \text{for each}\ i\in\set*{1,\dots,N},
\\[1ex]
\widehat u
&
\text{in}\
\widehat\Omega_r,
\\[1ex]
g
&
\text{in}\ 
\R^2\setminus\Omega.
\end{cases}
\end{equation}
As for $u_s$, we have that $\widetilde u_s\in H^1(\R^2;\R^2)$ and $\|\widetilde u_s\|_{L^\infty}\le\|g\|_{L^\infty}$.
In addition, by \cref{res:null_jac} (since $|\widetilde u_s|=1$ on $\widehat{(\R^2)}_r$) and by~\eqref{eq:vortex_jac} in \cref{res:vortex1}, we have 
\begin{equation*}
\jac(\widetilde u_s)
=
\sum_{i=1}^N
\jac(\upsilon_{d_i,0,(M_s+1)\tau_s}(\,\cdot\,-x_i))\,\chi_{B_{(M_s+1)\tau_s}(x_i)}
=
\pi
\sum_{i=1}^N
d_i
\,\frac{\chi_{B_{(M_s+1)\tau_s}(x_i)}}{|B_{(M_s+1)\tau_s}(x_i)|},
\end{equation*}
from which, recalling~\eqref{eq:special}, we deduce that
\begin{equation*}
\jac(\widetilde u_s)\xflute{A}\pi\mu
\quad
\text{as}\
s\to1^-.
\end{equation*} 
We can hence achieve~\eqref{eq:recovery_I_jac} by proving that 
\begin{equation}
\label{eq:recovery_I_aux}
\lim_{s\to1^-}
\|\jac(v_s)-\jac(\widetilde u_s)\|_{\flute(A)}
=0.
\end{equation}
To this end, we exploit \cref{res:flat_dist_jac}.
On the one hand, by~\eqref{eq:E_comp_L2} in \cref{res:E_comp} and~\eqref{eq:special_limsup}, and by~\eqref{eq:aux}, \eqref{eq:vortex_L2_dist} in \cref{res:vortex1}, and~\eqref{eq:tau_M}, we can estimate  
\begin{equation}
\label{eq:aux_j1}
\|v_s-\widetilde u_s\|_{L^2(A)}
\le 
\|v_s-u_s\|_{L^2(A)}
+
\|u_s-\widetilde u_s\|_{L^2(A)}
\le 
C\sqrt{(1-s)\,|\log(1-s)|}, 
\end{equation}
where $C>0$ is a constant which does not depend on~$s$.
On the other hand, by~\eqref{eq:lim_gamma_R}, \eqref{gradok}  and~\eqref{eq:special_limsup}, and by~\eqref{eq:aux} and~\eqref{eq:vortex_H1}, we also have that
\begin{equation}
\label{eq:aux_j2}
\|D v_s\|_{L^2(A)}+\|D \widetilde u_s\|_{L^2(A)}
\le
C
\sqrt{|\log(1-s)|},
\end{equation}
where $C>0$ is a constant which does not depend on~$s$.
Hence, by \cref{res:flat_dist_jac}, the validity of~\eqref{eq:recovery_I_aux} follows by combining~\eqref{eq:aux_j1} and~\eqref{eq:aux_j2}.

\vspace{1ex}

At this point, thanks to Steps~1 and~2 above, we proved~\ref{item:main_limsup} in \cref{resi:main} for any $\mu\in\atoms(\overline\Omega)$ as in~\eqref{eq:special} such that $\mu(\overline\Omega)=\deg(g|_{\partial\Omega},\partial\Omega)$.

\vspace{1ex}

\textit{Step~3}.
We now prove~\ref{item:main_limsup} in \cref{resi:main} in full generality.
Let $\mu\in\atoms(\overline\Omega)$ be such that $\mu(\overline\Omega)=\deg(g|_{\partial\Omega},\partial\Omega)$.
We can find a sequence $(\mu_k)_{k\in\N}\subset\atoms(\overline\Omega)$ as in~\eqref{eq:special}, that is, 
\begin{equation*}
\mu_k(\overline\Omega)=\deg(g|_{\partial\Omega},\partial\Omega)
\quad
\text{and}
\quad
\mu_k=\sum_{i=1}^{N_k}d_{i,k}\,\delta_{x_{i,k}}
\quad
\text{for each}\
k\in\N,
\end{equation*}
with $N_k\in\N$, $x_{i,k}\in\Omega$ and $d_{i,k}\in\set*{-1,1}$ for each $i\in\set*{1,\dots,N_k}$, such that 
\begin{equation*}
\mu_k\xflute{\overline\Omega}\mu 
\quad
\text{as}\ 
k\to\infty.
\end{equation*}
In particular, we have that $|\mu_k|(\Omega)\to|\mu|(\overline\Omega)$ as $k\to\infty$.
By Step~1, for each $k\in\N$ we can find $u_{s,k}\in H^s_g(\Omega;B_L)$ such that
\begin{equation*}
\limsup_{s\to1^-}
\GL_\lambda^s(u_{s,k})
\le 
\pi|\mu_k|(\Omega)
\end{equation*}
whenever $\lambda>0$.
Thus, given any sequence $(s_j)_{j\in\N}\subset(0,1)$ such that $s_j\to1^-$ as $j\to\infty$, by a diagonal argument we can extract a subsequence $(s_{j_k})_{k\in\N}$ such that 
\begin{equation*}
\lim_{k\to\infty}
\GL_\lambda^{s_{j_k}}(u_{s_{j_k},k})
\le 
\lim_{k\to\infty}
\pi|\mu_k|(\Omega)
=
\pi|\mu|(\overline\Omega),
\end{equation*} 
concluding the proof.
\end{proof}

%%% BIBLIO %%%%

\end{document}